\documentclass[11pt]{article}
\usepackage{times}
\usepackage{amsfonts,amscd,amssymb}
\newtheorem{lemma}{Lemma}[section]
\newtheorem{theorem}[lemma]{Theorem}
\newtheorem{proposition}[lemma]{Proposition}
\newtheorem{corollary}[lemma]{Corollary}
\newtheorem{definition}[lemma]{Definition}
\newtheorem{remark}[lemma]{Remark}
\newtheorem{remarks}[lemma]{Remarks}

\newtheorem{examples}[lemma]{Examples}

\newenvironment{proof}{{\it Proof.}}{\hfill $ \square $ \vskip 4mm}

\newcommand{\nc}{\newcommand}
\nc{\rnc}{\renewcommand}
\nc{\nt}{\newtheorem}

        {~\hfill~\fbox{}\end{trivlist}}

\nc{\thlabel}[1]{\label{th:#1}}
\nc{\thref}[1]{Theorem~\ref{th:#1}}
\nc{\selabel}[1]{\label{se:#1}}
\nc{\seref}[1]{Section~\ref{se:#1}}
\nc{\lelabel}[1]{\label{le:#1}}
\nc{\leref}[1]{Lemma~\ref{le:#1}}
\nc{\prlabel}[1]{\label{pr:#1}}
\nc{\prref}[1]{Proposition~\ref{pr:#1}}
\nc{\colabel}[1]{\label{co:#1}}
\nc{\coref}[1]{Corollary~\ref{co:#1}}
\nc{\relabel}[1]{\label{re:#1}}
\nc{\reref}[1]{Remark~\ref{re:#1}}
\nc{\exlabel}[1]{\label{ex:#1}}
\nc{\exref}[1]{Example~\ref{ex:#1}}
\nc{\delabel}[1]{\label{de:#1}}
\nc{\deref}[1]{Definition~\ref{de:#1}}
\nc{\eqlabel}[1]{\label{eq:#1}}
\nc{\eqref}[1]{(\ref{eq:#1})}
\nc{\csm}{\mbox{$\triangleright\!\!\!<$}}
\nc{\smc}{\mbox{$>\!\!\!\triangleleft$}}
\nc{\trr}{\triangleright}
\nc{\Hom}{\rm{Hom}\,}
\nc{\Mor}{\rm{Mor}\,}
\nc{\Aut}{\rm{Aut}\,}
\nc{\Ann}{\rm{Ann}\,}
\nc{\Ker}{\rm{Ker}\,}
\nc{\Trace}{\rm{Trace}\,}
\nc{\Char}{\rm{Char}\,}
\nc{\Mod}{\rm{Mod}\,}
\nc{\End}{\rm{End}\,}
\nc{\Spec}{\rm{Spec}\,}
\nc{\Span}{\rm{Span}\,}
\nc{\sgn}{\rm{sgn}\,}
\nc{\Id}{\rm{Id}\,}
\nc{\Com}{\rm{Com}\,}

\def\lan{\langle}
\def\ran{\rangle}
\def\eps{\varepsilon}

\def\nnat{\mathbb{ N}}

\def\nint{\mathbb{ Z}}

\nc{\dht}{\mbox{$\rightharpoonup\hspace{-2ex}\rightharpoonup$}}
\nc{\dhtb}{\mbox{$\leftharpoonup\hspace{-2ex}\leftharpoonup$}}
\nc{\nd}{\mbox{$\not|$}} %not divide sign
\def\text#1{\mbox{{\rm #1}}}

\nc{\nci}{\mbox{$\not\subseteq$}}
\nc{\scontainin}{\mbox{$\mbox{}\subseteq\hspace{-1.5ex}\raisebox{-.5ex}
{$_\prime $}\hspace*{1.5ex}$}}

\def\ot{\otimes}

\def\doublerightleft#1#2{{\lower.2ex\vbox{
\hbox{${\smash{\mathop{\longrightarrow}\limits^{#1}}}$}\vspace*{-4mm}
\hbox{${\smash{\mathop{\longleftarrow}\limits_{#2}}}$}}}}

\input diagrams

\begin{document}
\title{The structure of Frobenius algebras and separable algebras}
\author{S. Caenepeel \\ Faculty of Applied Sciences\\
Free University of Brussels, VUB
\\ Pleinlaan 2 \\ B-1050 Brussels, Belgium
\and
Bogdan Ion\\
Department of Mathematics\\
Princeton University\\
Fine Hall, Washington Road\\
Princeton, NJ 08544-1000, USA
\and
G. Militaru\thanks{Research supported by the bilateral project
``Hopf algebras and co-Galois theory" of the Flemish and
Romanian governments}
\\ Faculty of Mathematics\\ University of Bucharest\\ Str. Academiei 14
\\RO-70109 Bucharest 1, Romania }
\date{}
\maketitle
\begin{abstract}
\noindent We present a unified approach to the study of separable
and Frobenius algebras. The crucial observation is that both types
of algebras are related to the nonlinear equation
$R^{12}R^{23}=R^{23}R^{13}=R^{13}R^{12}$, called the FS-equation.
Solutions of the FS-equation automatically satisfy the braid equation,
an equation that is in a sense equivalent to the quantum
Yang-Baxter equation.
Given a solution to the FS-equation satisfying a certain normalizing
condition, we can construct a Frobenius algebra or a separable
algebra ${\cal A}(R)$ - the normalizing condition is different in
both cases. The main result of this paper is
the structure of these two fundamental types of algebras: a
finitely generated projective Frobenius or separable $k$-algebra $A$
is isomorphic to such an ${\cal A}(R)$. If $A$ is a free $k$-algebra,
then ${\cal A}(R)$ can be described using generators and relations.
A new characterization of Frobenius extensions is given:
$B\subset A$ is Frobenius if and only if $A$ has a $B$-coring
structure $(A, \Delta, \varepsilon)$ such that the comultiplication
$\Delta:\ A\to A\ot_B A$ is an $A$-bimodule map.
\end{abstract}

\section{Introduction}
Frobenius extensions in noncommutative ring theory have been
introduced by Kasch \cite{FK}, as generalizations of the classical
notion of Frobenius algebras over a field (see also \cite{N1},
\cite{N2}, \cite{PAR}).
The notion of Frobenius extension has turned out to be a
fundamental one, and many generalizations have appeared in the
literature. Roughly stated, an object ${\cal O}$ satisfies a
``Frobenius-type" property if two conditions hold:
a ``finiteness" condition (for example finite dimensionality) and a ``symmetry"
condition (for example
${\cal O}$ has an isomorphic dual ${\cal O}^*$).
${\cal O}$ can be an algebra (\cite{CR}), a coalgebra
(\cite{D}), a Lie algebra (\cite{Dr}), a Hopf algebra (\cite{P})
or, in the most general case, a functor between two categories
(\cite{CMZ}). In fact an object ${\cal O}$ over a field $k$ is
Frobenius if the forgetful functor
$F:\ {}_{\cal O}{\rm Rep} \to {}_k{\cal M}$ from the category of
representations of ${\cal O}$ to vector spaces is Frobenius, and this
means that $F$ has a left adjoint which is also a right adjoint.\\
In recent years, many important new results about ``Frobenius
objects" came about
(see \cite{A1}-\cite{A}, \cite{BDGN}-\cite{BF},
\cite{B1}, \cite{CMZ}, \cite{F}, \cite{FMS}, \cite{K1}). Let us mention
a few of them, illustrating the
importance of the concept. The homology of a compact oriented
manifold is a Frobenius algebra. In \cite{A1} (see also \cite{A}), is is
shown that
there is a one-to-one correspondence between two-dimensional
topological quantum field theories and Frobenius algebras.
Frobenius extensions also appear in quantum cohomology,
where they provide a generalization of the classical Euler class
(\cite{A2}).
In \cite{K1} the Jones polynomial is constructed starting
from certain separable Frobenius ring extensions, instead of starting
from type ${\rm II_1}$ subfactors as in the original work of
Jones.
In \cite{Dr} Lie quasi-Frobenius algebras are studied in connection to
the classical Yang-Baxter equation, while Frobenius algebras appear
in the study of the quantum Yang-Baxter equation (\cite{BFS}).
In \cite{Dub} the concept of Frobenius manifold is introduced and
discussed: a manifold $M$ is Frobenius if each fiber of the tangent
bundle $TM$ has a Frobenius algebra structure.
Frobenius functors were introduced in \cite{CMZ}. The definition is
based on the long-standing observation that a ring
extension $B\subseteq A$ is Frobenius if and only
if the restriction of scalars functor
${}_A{\cal M}\to {}_B{\cal M}$ has isomorphic left and right
adjoint functors (\cite{MN}, \cite{Mo}).
Theorem 4.2 in \cite{CMZ} states that the forgetful functor from
the category of Yetter-Drinfel'd modules over a Hopf algebra $H$ to
the category of $H$-modules is Frobenius if and only if
$H$ is finite dimensional and unimodular. This can be viewed as
the ``quantum version" of the classical result that any
finite dimensional Hopf algebra is Frobenius. We remark that, in
this situation, the symmetry condition that we mentioned earlier
comes down to unimodularity of $H$.\\
Another fundamental concept in the theory of (non)commutative rings
is that of separable extensions, we refer to \cite{DI} for a detailed
discussion. In \cite{NOB}, the notion of separable functor is
introduced, and, as we have explained in \cite{CIMZ}, the existing
versions of Maschke's Theorem come down to proving that a certain
functor is separable: an object ${\cal O}$ is semisimple (or reducible)
if and only if the forgetful functor
$F:\ {}_{\cal O}{\rm Rep} \to {}_k{\cal M}$, is separable.
For example, we can show that the above mentioned forgetful functor
from Yetter-Drinfel'd modules to $H$-modules is separable
if and only if there exists a so-called quantum integral, this is
a map with certain properties (see \cite[Theorem 3.13]{CIMZ}).
In \cite{B2}, separable functors are applied to the category of
representations of an entwined structure, introduced in \cite{BM}.\\
In this paper we study Frobenius algebras and separable
algebras from a unifying point of view. In order to do this,
we focus attention to what we have called the
FS-equation (or Frobenius-separability equation), namely
\begin{equation}\eqlabel{0.1}
R^{12}R^{23}=R^{23}R^{13}=R^{13}R^{12}
\end{equation}
where $A$ is an algebra over a commutative ring $k$ and
$R\in A\ot A$. Solutions of \eqref{0.1} automatically satisfy the
{\it braid equation}: $R^{12}R^{23}R^{12}=R^{23}R^{12}R^{23}$. Also observe
that,
in the case where $A=\End_k(M)$, $R$ is a solution of the braid
equation if and only if $\tau\circ R$ (or, equivalently, $R\circ\tau$)
is a solution of the {\it quantum Yang-Baxter equation}
$R^{12}R^{13}R^{23}=R^{23}R^{13}R^{12}$.\\
Let $A$ be an algebra with unit, and $\Delta:\ A\to A\ot A$
an $A$-bimodule map. Then $R=\Delta(1)$ is a solution of \eqref{0.1}.
It is well-known that $A$ is separable if such a $\Delta$ exists,
with a certain normalizing condition. Moreover, as shown recently
in \cite{A3}, $A$ is Frobenius if there exists a map $\Delta$ as
above, together with a counit map $\eps:\ A\to k$, satisfying a
different normalizing condition. So both the Frobenius property and
separability are related to solutions of the FS-equation, and this
explains our terminology. We will call the two different normalizing
conditions the normalizing separability condition and the
normalizing Frobenius condition. Given a separable or Frobenius
algebra, we find a solution of the FS-equation, satisfying one of
the two normalizing conditions.\\
The organisation of this paper is the following. In \seref{2},
we generalize a result of Abrams (see \cite[Theorem 2.1]{A3})
stating that a finite dimensional algebra $A$ over a field $k$
is Frobenius if and only if there exists a coalgebra structure on
$A$ such that the comultiplication is a bimodule map.
Based on the concept of Frobenius functors and drawing on the
Tannaka-Krein duality theory (see \cite{JS}, \cite{Maj 1}),
we can generalize this result to arbitrary ring extensions
(see \thref{2.4}). In \prref{2.6}, we give a criterion for
the separability of Frobenius extensions. Using duality
arguments, we obtain similar results for coalgebras.\\
In \seref{3}, we present a survey of solutions of the FS-equation.
Since solutions of the FS-equation are also solutions of the braid
equation, they can be transformed into solutions of the
quantum Yang-Baxter equation, leading to
interesting two-dimensional solvable integrable models in statistical
mechanics and quantum field theory (cf. \cite{Ba1}, \cite{Fa},
\cite{Hi}). We will give three new classes of solutions:
in \prref{3.7}, we construct solutions starting from a Sweedler
two-cocycle $\sigma:\ H\ot H\to k$, and we give explicit examples
in the case where $H$ is a finite group ring. In \prref{lazio}
we construct a solution for the $FS$-equation
starting from a map $\theta :\{1,\cdots, n\}^3\to \{1,\cdots, n\}$
satisfying certain properties. In \prref{3.14}
a solution of the $FS$-equation is associated to any idempotent
function $\phi: \{1, \cdots, n \} \to \{1, \cdots, n \}$.\\
\seref{4} contains our main result: given a
solution of the FS-equation, we can construct an algebra
${\cal A}(R)$ and a bimodule map $\Delta$ such that $R=\Delta(1)$.
If $R$ satisfies one of the two normalizing conditions, then
${\cal A}(R)$ is separable, resp. Frobenius (see \thref{4.1}).
Moreover, a finitely generated projective
$k$-algebra $A$ is separable or Frobenius, if it is isomorphic to
an algebra ${\cal A}(R)$, with $R$ a solution of the FS-equation
satisfying the appropriate normalizing condition (see \thref{4.2}).
In the case where $A$ is free of finite rank, ${\cal A}(R)$
can be described using generators and relations. Duality arguments
yield similar results about the structure of co-Frobenius and
coseparable coalgebras.\\
In \seref{5}, we introduce the category ${}_A{\cal FS}^A$ of
so-called $FS$-objects. These $FS$-objects can be compared to
Yetter-Drinfel'd modules in the situation where one works with the
braid equation rather than the $FS$-equation.
An interesting aspect of this category is that it measures in a
certain sense how far an algebra that is also a coalgebra is
different from a bialgebra: over a bialgebra, the category of
$FS$-objects happens to be just the category
of $k$-modules. Another property of the category ${}_A{\cal FS}^A$ is
that it is isomorphic with the category of left $A$-modules when $A$ has
a unit, respectively with the category of right $A$-comodules when $A$
has a counit. As an immediate consequence we obtain
\cite[Theorem 3.3]{A3}: for a Frobenius algebra $A$, the categories
${}_A{\cal M}$ and  ${\cal M}^A$ are isomorphic.

\section{Preliminaries}\selabel{1}
Let $k$ be a commutative ring. All modules,
algebras or coalgebras
are assumed to be  over $k$ and $\ot =\ot_k$.
For the comultiplication $\Delta$ on a $k$-coalgebra
$C$, we will use Sweedler's notation, that is
$$\Delta(c)=\sum c_{(1)}\ot c_{(2)}\in C\ot C$$
for any $c\in C$. If $M$ is a right $C$-comodule, then we will say that
$C$ coacts on $M$, and use the notation
$$\rho_M(m)=\sum m_{<0>}\ot m_{<1>}\in M\ot C.$$
${\cal M}^C$ will be the category of right $C$-comodules and $C$-colinear
maps. Let $(A, m_A, 1_A)$ be a $k$-algebra with associative
multiplication $m_A:\ A\ot A\to A$ and unit $1_A$. All the algebras
that we will consider in Sections 1-4 are associative with unit.
${}_A{\cal M}$ will be the category of
left $A$-modules and $A$-linear maps over the $k$-algebra $A$.

\subsection{Separable algebras and Frobenius algebras}\selabel{1.1}
For a $k$-module $M$, $I_M:\ M\to M$ and $\tau_M:\ M\ot M\to M\ot M$
will denote respectively the identity map on $M$, and the switch map,
mapping $m\ot n$ to $n\ot m$.\\
For a linear map $R:\ M\ot M\to M\ot M$, we will consider the maps
$$R^{12},~R^{13},~R^{23}:\ M\ot M\ot M\to M\ot M\ot M$$
defined by the formulas
$$R^{12}=R\ot I,~~~ R^{23}=I\ot R,~~~
R^{13}=(I\ot \tau)(R\ot I)(I\ot \tau)$$
A similar notation will be used for elements $e\in A\ot A$, where $A$ is a
$k$-algebra: if $e=\sum e^1\ot e^2$, then
$$e^{12}=\sum e^1\ot e^2\ot 1~~;~~
e^{13}=\sum e^1\ot 1\ot  e^2~~;~~
e^{23}=\sum 1\ot e^1\ot e^2$$
A sum over an empty family will be assumed to be zero:
$$\sum_{i\in \emptyset} a_i=0.$$
Let $(A, m_A, 1_A)$ be a $k$-algebra with unit. Then $A\ot A$ is an
$A$-bimodule with the natural actions
$$a\cdot (b\ot c)\cdot d:=ab\ot cd$$
for all $a$, $b$, $c$, $d\in A$. An element $e=\sum e^1\ot e^2\in A\ot A$
will be called $A$-{\sl central} if for any $a\in A$ we have
\begin{equation}\eqlabel{2.3.1}
a\cdot e= e\cdot a
\end{equation}
Let $B\subset A$ be an arbitrary ring extension.
Recall from \cite{DI} that $A/B$ is called separable if there exists a
separability idempotent, that is an $A$-central element
$e=\sum e^1\ot_B e^2 \in A\ot_B A$ satisfying the
normalizing separability condition
\begin{equation}\eqlabel{SNa}
\sum e^1e^2 = 1
\end{equation}

\begin{remark}\relabel{1.1}
\rm If $B\subset A$ is a separable extension, then $A/B$ is a left
(right) semisimple extension in the sense of \cite{HS}: this means that
any
left $A$-submodule $W$ of a left $A$-module $V$ having a
$B$-complement in $V$, has an $A$-complement in $V$. Indeed, let
$e=\sum e^1\ot_B e^2 \in A\ot_B A$ be a separability idempotent and
$f:\ V\to W$ be a left $B$-module map splitting the inclusion
$W\subset V$. Then
$$\tilde{f}:\ V\to W, \quad \tilde{f}(v):=\sum e^1f(e^2 v)$$
for all $v\in V$, is a left $A$-module map splitting the inclusion
$W\subset V$.\\
Observe that the deformation $f\to \tilde{f}$ is functorial. As we
have explained in \cite{CIMZ}, this is what happens in almost all
generalizations of Maschke's Theorem that have appeared in the
literature: one constructs a deformation of a map, and this
deformation turns out to be functorial.
\end{remark}

Let $B\subset A$ be a ring extension. Then
$\Hom_B(A,\; B)$ is an $(A, B)$-bimodule with left $A$-action and
right $B$-action given by the formulas
$$(a\cdot f \cdot b)(x)= f(xa)b$$
for all $a$, $x\in A$, $b\in B$ and $f\in \Hom_B(A,\; B)$.
>From \cite{N1}, we recall the following

\begin{definition}\delabel{1.2}
Let $B\subset A$ be a ring extension. Then $A/B$ is called
a Frobenius extension (or, a Frobenius extension of the first kind) if
the following two conditions hold:\\
\hspace*{7mm}1) ${}_BA$ is finitely generated and projective in ${}_B{\cal
M}$;\\
\hspace*{7mm}2) there exists an isomorphism
$$\varphi:\ A \to \Hom_B (A,\; B)$$
\hspace*{7mm}of $(A,B)$-bimodules.
\end{definition}

The concept of Frobenius extension is left-right symmetric: that is
any left Frobenius extension is also a right Frobenius extension
(see \cite[Proposition 1]{N1}). Of course, for $B=k$ we obtain the
clasical concept of Frobenius algebra: a $k$-algebra $A$ is Frobenius
if $A$ is finitely generated and projective as a $k$-module, and
$A\cong A^*$ as left $A$-modules.\\
A $k$-coalgebra $C$ is called coseparable (\cite{L}) if there exists a
coseparability idempotent, that is a $k$-linear map
$\sigma:\ C\ot C\to k$ such that
$$\sum \sigma(c\ot d_{(1)})d_{(2)}=
\sum \sigma(c_{(2)}\ot d)c_{(1)}~~{\rm and}~~
\sum \sigma(c_{(1)}\ot c_{(2)})=\varepsilon(c)$$
for all $c$, $d\in C$.

\subsection{Matrix algebras and comatrix coalgebras}\selabel{1.2}
Let $M$ be a free $k$-module, with basis $\{m_1,\cdots,m_n\}$,
and $\{p^1,\cdots,p^n\}$ the corresponding dual basis of the
dual module $M^*$, such that
\begin{equation}\eqlabel{matrix1}
\lan p^i,m_j\ran =\delta^i_j
\end{equation}
for all $i,j\in\{1,\cdots,n\}$. $\delta^i_j$ is the Kronecker symbol.
Then $\{e^i_j=p^i\ot m_j~|~i,j=1,\cdots,n\}$ and
$\{c^i_j=m_j\ot p^i~|~i,j=1,\cdots,n\}$ are free basis for
respectively $\End_k(M)\cong M^*\ot M$ and
$\End_k(M^*)\cong M\ot M^*$. The isomorphisms are given by the rules
\begin{equation}\eqlabel{matrix2}
e^i_j(m_k)=\delta^i_km_j~~{\rm and}~~c^i_j(p^k)=\delta^k_j p^i
\end{equation}
$\End_k(M)\cong M^*\ot M$ is isomorphic to the $n\times n$-matrix algebra
${\cal M}_n(k)$, and the multiplication in $M^*\ot M$ is given by the rules
\begin{equation}\eqlabel{matrix3}
e^i_je^k_l=\delta^i_le^k_j~~{\rm and}~~\sum_{i=1}^n e^i_i=1
\end{equation}
$\End_k(M^*)\cong M\ot M^*$ is isomorphic to the $n\times n$-comatrix
coalgebra ${\cal M}^n(k)$, and the comultiplication on $M\ot M^*$
is given by
\begin{equation}\eqlabel{matrix4}
\Delta(c^i_j)=\sum_k c^i_k\ot c^k_j~~{\rm and}~~\varepsilon(c^i_j)=\delta^i_j
\end{equation}
We obtain the matrix algebra ${\cal M}_n(k)$ and the comatrix coalgebra
${\cal M}^n(k)$ after we take $M=k^n$ and the canonical basis.\\
A linear map $R:\ M\ot M\to M\ot M$ can be described by a matrix $X$,
with $n^4$ entries $x^{ij}_{uv}\in k$, where $i,j,u,v$ range from
$1$ to $n$. This means
\begin{equation}\eqlabel{matrix5}
R(m_u\ot m_v)=\sum_{ij} x^{ij}_{uv} m_i\ot m_j
\end{equation}
or
\begin{equation}\eqlabel{matrix6}
R=\sum_{ijuv} x^{ij}_{uv} e^u_i\ot e^v_j
\end{equation}

\section{Frobenius functors}\selabel{Ff}\selabel{2}
In \cite[Theorem 2.1]{A3}, L. Abrams gives an interesting
characterization of Frobenius algebras. Let $k$ be a field
and $A$ a finite dimensional $k$-algebra. Then $A$ is Frobenius
if and only if there exists a coalgebra structure $(A, \Delta, \eps)$
on $A$ such that the comultiplication $\Delta$ is an $A$-bimodule map.
The coalgebra structure which arises on a Frobenius algebra
can be briefly described as follows:
let $A$ be a Frobenius algebra and $\varphi:\ A \to A^*$ be a left
$A$-module isomorphism. Let $(A^*, m_{A^*}, 1_{A^*})$ be the unique
algebra structure on $A^*$ such that $\varphi$ becomes an isomorphism
of algebras. Using the usual duality between the categories of finite
dimensional
algebras and finite dimensional coalgebras
(see \cite{Ab}), we obtain that $A$ has a unique coalgebra structure
such that ${\cal M}^A \cong {}_{A^*}{\cal M}\cong {}_{A}{\cal M}$.
In particular, \cite[Theorem 3.3]{A3} follows.
In this Section, we will generalize this result to arbitrary ring
extensions $B\subset A$. The methods of \cite{A3} do not work in this
case; our approach is based on the notion of Frobenius functor, as
introduced in \cite{CMZ}. Our construction of the comultiplication $\Delta$
and the
counit $\eps$ associated to a Frobenius extension is inspired by the
Tannaka-Krein
duality theory (see \cite{JS}, \cite{Maj 1}). We begin with the following
key result.

\begin{proposition}\prlabel{nori}\prlabel{2.1}
Let $(A,m_A,1_A)$ be a $k$-algebra, $\Delta:\ A\to A\ot A$ an
$A$-bimodule map and $e=\Delta(1_A)$.\\
1) In $A\ot A\ot A$, we have the equality
\begin{equation}\eqlabel{2.3.2a}
e^{12}e^{23}=e^{23}e^{13}=e^{13}e^{12}
\end{equation}
2) $\Delta$ is coassociative;\\
3) If $(A, \Delta, \eps)$ is a coalgebra structure on $A$, then
$A$ is finitely generated and projective as a $k$-module.
\end{proposition}

\begin{proof}
1) From the fact that $\Delta$ is an $A$-bimodule map,
it follows immediately that
$e$ $A$-central.
Write $E=\sum E^1\ot E^2=e$. Then
\begin{eqnarray*}
e^{12}e^{23}&=& \sum (e^1\ot e^2\ot 1)(1\ot E^1\ot E^2)\\
&=&\sum \underline{e^1\ot e^2E^1 }\ot E^2\\
\text{($\mbox{using} \;\eqref{2.3.1}$)}
&=&\sum E^1e^1\ot e^2\ot E^2=e^{13}e^{12}\\
{\rm and}~~~~~~~~~~~~~~~~&&\\
e^{12}e^{23}&=& \sum e^1\ot \underline{e^2E^1 \ot E^2}\\
\text{($\mbox{using} \; \eqref{2.3.1}$)}
&=&\sum e^1\ot E^1\ot E^2e^2=e^{23}e^{13}.
\end{eqnarray*}
2) follows from 1) and the formulas
\begin{eqnarray*}
(\Delta\ot I)\Delta(a)&=&e^{12}e^{23}\cdot (1_A\ot 1_A\ot a)\\
(I\ot \Delta)\Delta(a)&=&e^{23}e^{13}\cdot (1_A\ot 1_A\ot a)
\end{eqnarray*}
for all $a\in A$.

3) Let $a\in A$. Applying $\eps\ot I_A$ and $I_A\ot\eps$ to
\eqref{2.3.1}, we obtain, using the fact that $\eps$ is a counit map,
$$a=\sum \eps (ae^1)e^2=\sum e^1\eps (e^2a)$$
and it follows that $\{e^1, \eps(e^2 \bullet)\}$
(or  $\{e^2, \eps(\bullet e^1)\}$) are dual bases of $A$
as a $k$-module.
\end{proof}

Let ${\cal C}$ and ${\cal D}$ be two categories, and let
$$\begin{diagram}
{\cal C} & \pile{ \rTo^F \\ \lTo_G} & {\cal D},
\end{diagram}$$
be two covariant functors. Observe that we have four covariant functors
$$\Hom_{\cal C}(G(\bullet),\bullet):\
{\cal D}^{\rm op}\times {\cal C}\to \underline{\rm Sets}
~~,~~
\Hom_{\cal D}(\bullet,F(\bullet)):\ {\cal D}^{\rm op}\times {\cal C}\to
\underline{\rm Sets},
$$
$$
\Hom_{\cal C}(\bullet,G(\bullet)):\
{\cal C}^{\rm op}\times {\cal D}\to \underline{\rm Sets}
~~{\rm and}~~
\Hom_{\cal D}(F(\bullet),\bullet):\ {\cal C}^{\rm op}\times {\cal D}\to
\underline{\rm Sets}.
$$
Recall that $G$ is a left adjoint for $F$ if
$\Hom_{\cal C}(G(\bullet),\bullet)$ and
$\Hom_{\cal D}(\bullet,F(\bullet))$ are isomorphic. $G$ is a right
adjoint for $F$ if the functors
$\Hom_{\cal C}(\bullet,G(\bullet))$ and
$\Hom_{\cal D}(F(\bullet),\bullet)$
are isomorphic.

\begin{definition} (\cite{CMZ})\delabel{2.2}
A covariant functor $F:\ {\cal C}\to {\cal D}$ is called a Frobenius functor if
it has isomorphic left and right adjoints.
\end{definition}

\begin{remark}\relabel{bel}\relabel{2.3}
\rm
Let us explain the terminology. Take an arbitrary
ring extension $B\subset A$ and let
$$F:\ {}_A{\cal M}\to {}_B{\cal M}$$
be the restriction of scalars functor. Then $F$ has a left adjoint,
the induction functor,
$$G=A\ot_B \bullet :\ {}_B{\cal M}\to {}_A{\cal M}$$
and a right adjoint, the coinduction functor,
$$T=\Hom_B(A,\bullet):\ {}_B{\cal M}\to {}_A{\cal M}$$
The left $A$-module structure on $T(N)=\Hom_B(A,N)$ is given
by $(a\cdot \gamma)(b)=\gamma(ba)$, for all $a,b\in A$ and
$\gamma\in \Hom_B(A,N)$. For a left $B$-module map $F$, we put
$T(f)=f\circ\bullet$.\\
The restriction of scalars functor
$F:\ {}_A{\cal M}\to {}_B{\cal M}$ is a Frobenius functor
if and only if the extension $B\subset A$ is Frobenius.
Indeed, the adjoint of a functor is unique upto
isomorphism (see \cite{Kan}). Hence, $F$ is a Frobenius functor
if and only if the induction functor $G$ is isomorphic to the
coinduction functor $T$.
This is equivalent (cf. \cite[Theorem 2.1]{Mo}, \cite[Theorem 3.15]{MN})
to the fact that $B\subset A$ is a Frobenius extension.
For a further study of Frobenius functors we refer to \cite{CGN}.
\end{remark}

In the next Theorem, $B\subset A$ will be an arbitrary ring extension.
The multiplication $m_A:\ A\times A\to A$ is then $B$-balanced, which
means that
it induces a map $A\ot_B A\to A$, which is also denoted by $m_A$.

\begin{theorem}\thlabel{abc}\thlabel{2.4}
Let $B\subset A$ be an arbitrary ring extension.
Then the following statemens are equivalent:
\begin{enumerate}
\item[1)] the restriction of scalars functor
$F:{}_A{\cal M}\to {}_B{\cal M}$
is Frobenius (or, equivalently, $B\subset A$ is a Frobenius extension);
\item[2)] there exists a pair $(\Delta, \eps)$ such that
\begin{enumerate}
\item $\Delta:\ A\to A\ot_B A$ is an $A$-bimodule map;
\item $\eps:\ A\to B$ is a $B$-bimodule map, and a counit for
$\Delta$, that is
$$(I_A\ot_B \eps)\Delta=(\eps \ot_B I_A)\Delta=I_A.$$
\end{enumerate}
\end{enumerate}
\end{theorem}

\begin{proof}
Let $G,~T:\ {}_B{\cal M}\to {}_A{\cal M}$ be as above.\\
$1) \Rightarrow 2).$ Assume
that $G$ and $T$ are isomorphic and let $\phi:\ G=A\ot_B \bullet \to
T=\Hom_B(A,
\bullet )$ be a natural transformation, defining an isomorphism of functors.
We then have the following isomorphisms of left $A$-modules:
$$\phi_{A}:\ A\ot_B A \to \Hom_B(A,A)~~~{\rm and}~~~
\phi_{B}:\ A(\cong A\ot_B B)\to \Hom_{B}(A,\; B)$$
For $a\in A$, we define $\varphi_a:\ A\to A$ by $\varphi_a(b)=ba$.
Then the map $\varphi:\ A\to \Hom_B(A, A)$ mapping $a$ to $\varphi_a$
is left $A$-linear. Now take
$$\Delta =\phi^{-1}_{A}\circ \varphi~({\rm
i.e.~}\Delta(a)=\phi_A^{-1}(\varphi_A))
~~~\mbox{and}~~~
\eps=\phi_{B}(1_A)$$
$\Delta$ is left $A$-linear, because $\phi^{-1}_{A}$ and $\varphi$
are left $A$-linear. $\phi$ is a natural transformation, so we
have a commutative diagram
$$\begin{diagram}
\Hom_B(A,A) & \rTo^{\phi^{-1}_A} & A\ot_B A \\
\dTo^{\varphi_a\circ\bullet} &  & \dTo_{I_A\ot \varphi_a}\\
\Hom_B(A,A) & \rTo_{\phi^{-1}_A} & A\ot_B A
\end{diagram}$$
for any $a\in A$. Applying the diagram to the identity map $I_A=\varphi_1$,
we find that $\Delta(a)=\Delta(1)a$, and $\Delta$ is also right
$A$-linear.\\
By construction, $\eps$ is left $B$-linear.
We will now prove that $\eps$ is also right $B$-linear,
and a counit for $\Delta$. Take
$N\in {}_B{\cal M}$ and $n\in N$, and consider the unique left
$B$-linear map $\varphi_n:\ B\to N$ mapping $1_B$ to $n$.
$\phi$ is a natural transformation, so we have a commutative diagram
$$\begin{diagram}
A\ot_B B & \rTo^{\phi_B} & \Hom_B(A,B) \\
\dTo^{I_A\ot\varphi_n} &  & \dTo_{\varphi_n\circ\bullet}\\
A\ot_B N & \rTo_{\phi_N} &  \Hom_B(A,N)
\end{diagram}$$
and it follows that $\phi_N(1_A\ot_B n)=\varphi_n\circ \eps$ and
$$\phi_N(1_A\ot_B n)(a)=\eps(a)n$$
Using the left $A$-linearity of $\phi_N$, we obtain
\begin{eqnarray}
\phi_N(a\ot_B n)(b) &=& \phi_N(a\cdot(1_A\ot_B n))(b)\nonumber\\
&=& (a\cdot \phi_N(1_A\ot_B n))(b)\nonumber\\
&=&\phi_N(1_A\ot_B n)(ba)\nonumber\\
&=& \eps(ba)n\eqlabel{gen}.
\end{eqnarray}
for all $a$, $b\in A$ and $n\in N$. In particular, for $N=B$, $b=1$ and
$n=\beta\in B$ we obtain that
$$\phi_B(a\beta\ot_B 1)(1)=\eps(a)\beta.$$
Hence,
$$\eps(a\beta)= (a\beta \cdot \eps)(1)=(a\beta \cdot \phi_B(1_A\ot 1_B))(1)=
\phi_B(a\beta\ot 1_B)(1)=\eps(a)\beta,$$
and $\eps$ is also right $B$-linear.

Applying \eqref{gen} with $a\ot_B n$ replaced by
$\Delta(1_A)=\sum e^1\ot_B e^2$, and $b$ by $1_A$, we find
$$\phi_A(\sum e^1\ot_B e^2)(1_A)=\sum \eps(e^1)e^2.$$
On the other hand,
$$\phi_A(\Delta(1_A))(1_A)=\varphi_1(1_A)=I_A(1_A)=1_A$$
and it follows that $\sum \eps(e^1)e^2=1_A$.\\
$\varepsilon:\ A\to B$ is left $B$-linear, and $\phi$ is a
natural transformation, so we have
a commutative diagram
$$\begin{diagram}
\Hom_B(A,A)& \rTo^{\phi^{-1}_A} & A\ot_B A  \\
\dTo^{\eps\circ\bullet} &  & \dTo_{I_A\ot_B\eps} \\
\Hom_B(A,B)& \rTo_{\phi^{-1}_B} & A\ot_B B
\end{diagram}$$
Applying the diagram to $I_A=\varphi_1$, we find that $\sum e^1\eps(e^2)=1_A$.
We conclude that
$$\sum \eps(e^1)e^2=\sum e^1\eps(e^2)=1_A$$
and, using the fact that $\Delta$ is an $A$-bimodule map, we see
that $\eps$ is a counit.

$2) \Rightarrow 1).$ Assume $\Delta:\ A\to A\ot_B A$ and
$\eps:\ A\to B$ satisfy the conditions of the second statement
in the Theorem, and write
$\Delta(1_A)=e=\sum e^1\ot_B e^2\in A\ot_B A$. $e$ is $A$-central
since $\Delta$ is an $A$-bimodule map. We will construct two
natural transformations
$$\begin{diagram}
G & \pile{ \rTo^{\phi} \\ \lTo_{\psi}} & T ,
\end{diagram}$$
and prove that they are each others inverses.
For $N\in {}_B{\cal M}$, we define
$$
\begin{diagram}
A\ot_B N & \pile{ \rTo^{\phi_{N}} \\ \lTo_{\psi_{N}}} & \Hom_{B}(A,N),
\end{diagram}
$$
$$
\phi_N(a\ot_B n)(b)=\eps(ba)n, ~~~~~~
\psi_N(\chi)=\sum e^1\ot_B \chi(e^2),
$$
for all $\chi\in \Hom_B(A,N)$, $a,b\in A$, and $n\in N$.
We leave it as an exercise to the reader to show that $\phi_N$ and $\psi_N$
are well-defined, and that $\phi$ and $\psi$ are natural transformations.
Let us show that $\phi$ and $\psi$ are each others inverses.
Applying
$m_A\circ (\eps \ot_B I_A)\circ (I_A\ot \chi)$
to the identity
$$\sum ae^1\ot_B e^2=\sum e^1\ot_B e^2a$$
we obtain, using the fact that $\eps$ is a counit,
\begin{equation}\eqlabel{347}
\sum \eps(ae^1)\chi (e^2)=\sum \eps(e^1)\chi (e^2a)=
\sum \chi(\eps (e^1)e^2a)=\chi(a)
\end{equation}
for all $a\in A$. Now
\begin{eqnarray*}
\phi_N(\psi_N(\chi))(a) &=& \sum \phi_N(e^1\ot_B \chi(e^2))(a) \\
&=& \sum \eps(ae^1) \chi(e^2)\\
\text{\eqref{347}}~~~~&=& \chi (a)
\end{eqnarray*}
and
\begin{eqnarray*}
\psi_N(\phi_N(a\ot_B n)) &=& \sum e^1\ot_B \psi_N(a\ot_B n)(e^2) \\
&=& \sum e^1\ot_B \eps(e^2a)n\\
\text{($e$ is $A$-central)}~~&=& \sum ae^1\ot_B \eps(e^2)n\\
\text{($\eps$ is a counit)}~~&=& a\ot_B n
\end{eqnarray*}
proving that $\phi$ and $\psi$ are each others inverses.
\end{proof}

\begin{remark}
\rm
Let $B$ be a ring. A $B$-{\sl coring} is (cf. \cite{K1}, \cite{W})
a triple $(U, \Delta, \eps)$, where $U$ is a $B$-bimodule,
$\Delta :U\to U\ot_B U$, $\eps :U\to B$ are $B$-bimodule maps such
that $\Delta$ is coassociative and $\eps$ is a counit for $\Delta$.
Using this concept, the above theorem can be restated as follows:
a ring extension $B\subset A$ is Frobenius if and only if $A$ has a
$B$-coring structure such that the comultiplication $\Delta$ is an
$A$-bimodule map.
\end{remark}

For $k$ a field,
the equivalence $1) \Leftrightarrow 2)$ in the next Corollary
has been proved recently by Abrams (see \cite[Theorem 2.1]{A3}).

\begin{corollary}\colabel{adi}
For an algebra $A$ over a commutative ring $k$,
the following statements are equivalent:
\begin{enumerate}
\item[1)] $A$ is a Frobenius algebra;
\item[2)] there exist a coalgebra structure
$(A, \Delta_A, \varepsilon_A)$ on $A$ such that the comultiplication
$\Delta_A:\ A\to A\ot A$ is an $A$-bimodule map;
\item[3)] there exist $e=\sum e^1\ot e^2 \in A\ot A$ and $\eps\in A^*$ such
that
$e$ is $A$-central and the normalizing Frobenius condition
\begin{equation}\eqlabel{FNa}
\sum \eps(e^1)e^2=\sum e^1\eps(e^2) =1_A
\end{equation}
is satisfied. $(e,\varepsilon)$ is called a Frobenius pair.
\end{enumerate}
\end{corollary}

\begin{proof}
$1) \Leftrightarrow 2)$ follows \thref{2.4}, if we
take into account that a bimodule map
$\Delta:\ A\to A\ot A$ is coassociative (\prref{nori}).\\
$2) \Leftrightarrow 3)$: observe that an $A$-bimodule map $\Delta:\ A\to
A\ot A$
is completely determined by $e=\Delta(1_A)$, and that there is
a bijective corespondence
between the set of all $A$-central elements and the set of
all $A$-bimodule maps $\Delta :A\to A\ot A$, and use the second
statement in \prref{nori}. It is easy to see that the counit property
is satisfied if and only if \eqref{FNa} holds.
\end{proof}

Let $A$ be a Frobenius algebra over a field $k$. In \cite{A2},
$\omega_A= (m_A \circ\Delta) (1_A) \in A$ is called
the {\it characteristic element} of $A$. $\omega_A$ generalizes
the classical Euler class $e(X)$ of a connected orientated compact
manifold $X$.
A commutative Frobenius algebra is semisimple if and only if
its characteristic element $\omega_A$ is invertible in $A$
(see \cite[Theorem 3.4]{A2}).\\
Let $B\subset A$ be an arbitrary Frobenius extension, and let $\Delta$ and
$\eps$ be as in \thref{abc}. We call
$\omega_{A/B}= (m_A \circ\Delta) (1_A) \in A$ the
{\it characteristic element} of the extension $B\subset A$.
Using this element, we give a criterion for the separability of
a Frobenius extension (recall from \reref{1.1} that any separable
algebra over a field $k$ is semisimple).

\begin{proposition}\prlabel{char}\prlabel{2.6}
Let $B\subset A$ be an arbitrary Frobenius extension such that
the characteristic element $\omega_{A/B}$ is a unit of $A$.
Then $B\subset A$ is a separable extension.
\end{proposition}

\begin{proof}
Let $\Delta(1_A)=\sum e^1\ot_B e^2$. Then
$$\sum ae^1\ot_B e^2=\sum e^1\ot_B e^2a$$
for all $a\in A$. Hence, $\omega_{A/B}\in {\rm Z}(A)$, the center of $A$.
It follows that its inverse $\omega_{A/B}^{-1}$ is also an element in
${\rm Z}(A)$. Now,
$$R:=\omega_{A/B}^{-1}\sum e^1\ot_B e^2$$
is a separability idempotent, i.e. $A/B$ is a separable extension.
\end{proof}

Let $C$ be a projective coalgebra. In
In \cite[Corollary 2.6]{CMZ} it
was proved that the forgetful functor $F:\ {}^C{\cal M}\to {}_k{\cal M}$ is
Frobenius
if and only if $C$ is finitely generated over $k$
and the dual algebra $C^*$ is a Frobenius algebra.
We will now give the coalgebra version of \coref{adi}.

\begin{theorem}\thlabel{abcd}\thlabel{2.7}
For a projective $k$-coalgebra $C$, the following statements are
equivalent:
\begin{enumerate}
\item[1)] the forgetful functor $F:{}^C{\cal M}\to {}_k{\cal M}$
is Frobenius;
\item[2)] $C$ is finitely generated and there exists an associative and
unitary algebra structure $(C, m_C, 1_C)$ on $C$ such that the
multiplication map $m_C:\ C\ot C\to C$ is a $C$-bicomodule map.
\end{enumerate}
\end{theorem}

\begin{proof}
If $F:\ {}^C{\cal M}\to {}_k{\cal M}$ is Frobenius then $C$ is finitely
generated (see \cite[Corollary 2.6]{CMZ}). The rest of the proof
follows from \thref{abc}, using the fact that the categories ${}^C{\cal M}$ and
${}_{C^*}{\cal M}$ are isomorphic if $C$ is a finitely
generated projective coalgebra (see \cite{Ab}).
\end{proof}

\begin{remarks}\relabel{2.7b}\rm
It seems natural to call a $k$-coalgebra $C$ co-Frobenius if and only if
the forgetful functor $F:{}^C{\cal M}\to {}_k{\cal M}$ is Frobenius.
Let us point out that a different definition appears in the literature
(see \cite{Lin}): $C$ is called right co-Frobenius if there exists a
right $C^*$-monomorphism from $C$ to $C^*$. If $C$ is projective, and
the forgetful functor is Frobenius, then $C$ is right c-Frobenius
(see \cite[Corollary 2.6]{CMZ}).\\
2) In \prref{3.10b}, we will present a necessary and sufficient condition
for the forgetful functor to be Frobenius, which also holds in the
situation where $C$ is not projective.
\end{remarks}

\section{The Frobenius-separability equation}\selabel{3}
\prref{nori} leads us to the following definition
\begin{definition}\delabel{3.1}
Let $(A, m_A, 1_A)$ be a $k$-algebra and $R=\sum R^1\ot R^2\in A\ot A$.\\
1) $R$ is called a solution of the FS-equation if
\begin{equation}\eqlabel{3.1.1}
R^{12}R^{23}=R^{23}R^{13}=R^{13}R^{12}
\end{equation}
in $A\ot A\ot A$.\\
2) $R$ is called a solution of the S-equation
if $R$ is a solution of the FS-equation and the
the normalizing separability separability condition holds:
\begin{equation}\eqlabel{SN}
\sum R^1R^2=1_A
\end{equation}
3) $(R,\eps)$ is called a solution of the F-equation
if $R$ is a solution of the $FS$-equation, and
$\eps \in A^*$ is such that
the normalizing Frobenius condition holds:
\begin{equation}\eqlabel{FN}
\sum \eps(R^1)R^2=\sum R^1\eps(R^2)=1_A
\end{equation}
4) Two solutions of the FS-equation are called equivalent if there
exists an invertible element $u\in A$ such that $S=(u\ot u)R(u^{-1}\ot
u^{-1})$.
\end{definition}

\begin{remarks}\relabel{3.2}\rm
1) The FS-equation appeared first in \cite[Lemma 3.6]{BFS}.
If $R$ is a solution of the FS-equation then $R$ is also
a solution of the braid equation
\begin{equation}\eqlabel{braid}
R^{12}R^{23}R^{12}=R^{23}R^{12}R^{23}
\end{equation}
2) In many applications, the algebra $A$ is of the form $A=\End_k(M)$,
where $M$ is a (projective) $k$-module. Then we can view $R$ as
an element of $\End_k(M\ot M)$, using the embedding
$\End_k(M)\ot \End_k(M)\subseteq \End_k(M\ot M)$. \eqref{3.1.1}
can then be viewed as an equation in $\End_k(M\ot M\ot M)$. This is
what we will do in \prref{3.3} and in most of the examples in this
Section. If $M$ is finitely generated and projective, then
we can identify $\End_k(M)\ot \End_k(M)$ and $\End_k(M\ot M)$.
In this situation, the braid equation \eqref{braid} can be transformed
into the quantum Yang-Baxter equation
\begin{equation}\eqlabel{QYB}
R^{12}R^{13}R^{23}=R^{23}R^{13}R^{12}
\end{equation}
More precisely, $R$ is a solution of \eqref{braid} if and only if
$\tau\circ R$ is a solution of \eqref{QYB} if and only if
$R\circ\tau$ is a solution of \eqref{QYB}
(see \cite[Remark 10.1.2]{Montgomery}).
\end{remarks}

The FS-equation can be rewritten as a matrix equation. Let $M$ be
a free module with basis $\{m_1,\cdots,m_n\}$. Using the notation
introduced in \seref{1.2}, we can represent a linear map
$R\in \End_k(M\ot M)$
by its matrix $X$, with entries $x^{ij}_{uv}\in k$,
with $i,j,u,v\in \{1,\cdots,n\}$ (see \eqref{matrix5} and \eqref{matrix6}).
It is then straightforward to compute
\begin{eqnarray}
R^{12}R^{23}(m_u\ot m_v\ot m_w)&=&
\sum_{i,j,l}\left(\sum_k x^{ij}_{uk}x^{kl}_{vw}\right) m_i\ot m_j\ot m_l
\eqlabel{3.3.0a}\\
R^{23}R^{13}(m_u\ot m_v\ot m_w)&=&
\sum_{i,j,l}\left(\sum_k x^{jl}_{vk}x^{ik}_{uw}\right) m_i\ot m_j\ot m_l
\eqlabel{3.3.0b}\\
R^{13}R^{12}(m_u\ot m_v\ot m_w)&=&
\sum_{i,j,l}\left(\sum_k x^{il}_{kw}x^{kj}_{uv}\right) m_i\ot m_j\ot m_l
\eqlabel{3.3.0c}
\end{eqnarray}

\begin{proposition}\prlabel{3.3}
Let $M$ be a free $k$-module with basis $\{m_1, \cdots, m_n \}$,
and let $R\in \End_k(M\ot M)$ be given by \eqref{matrix5}.\\
1) $R$ is a solution of the FS-equation if and only if
\begin{equation}\eqlabel{3.3.1}
\sum_k x^{ij}_{uk}x^{kl}_{vw}=
\sum_k x^{jl}_{vk}x^{ik}_{uw}=
\sum_k x^{il}_{kw}x^{kj}_{uv}
\end{equation}
for all $i,j,l,u,v,w\in\{1,\cdots,n\}$.\\
2) $R$ satisfies \eqref{SN} if and only if
\begin{equation}\eqlabel{3.3.2}
\sum_k x^{kj}_{ik} =\delta^j_i
\end{equation}
for all $i,j\in\{1,\cdots,n\}$.\\
3) Let $\varepsilon$ be the trace map. Then
$(R,\varepsilon)$ satisfies \eqref{FN} if and only if
\begin{equation}\eqlabel{3.3.3}
\sum_k x^{kj}_{ki} =\sum_k x^{jk}_{ik}=\delta^j_i
\end{equation}
for all $i,j\in\{1,\cdots,n\}$.
\end{proposition}

\begin{proof}
1) follows immediately from (\ref{eq:3.3.0a}-\ref{eq:3.3.0c}). 2) and 3)
follow from \eqref{matrix6}, using the multiplication rule
$$e^i_je^k_l=\delta^i_le^k_j$$
and the formula for the trace
$$\varepsilon(e^i_j)=\delta^i_j$$
\end{proof}

As a consequence of \prref{3.3}, we can show that
$R\in \End_k(M\ot M)$ is a solution of the equation
$R^{12}R^{23}=R^{13}R^{12}$ if and only if a certain multiplication on
$M\ot M$ is associative.

\begin{corollary}\colabel{3.4}
Let $M$ be a free $k$-module with basis $\{m_1, \cdots, m_n \}$,
and $R\in \End_k(M\ot M)$, given by \eqref{matrix5}.
Then $R$ is a solution of the equation $R^{12}R^{23}=R^{13}R^{12}$
if and only if the multiplication on $M\ot M$ given by
$$(m_k\ot m_l)\cdot (m_r\ot m_j):=
\sum_a x_{jl}^{ak}m_r\ot m_a$$
($k$, $l$, $r$, $j=1,\cdots,n$) is associative.\\
In this case $M$ is a left $M\ot M$-module with
structure
$$(m_k\ot m_l)\bullet m_j:=\sum_a x_{jl}^{ak}m_a$$
for all $k$, $l$, $j=1,\cdots,n$.
\end{corollary}

\begin{proof}
Write $m_{kl}=m_k\ot m_l$. Then
$$m_{pq}\cdot (m_{kl}\cdot m_{rj})=\sum_a x_{jl}^{ak} m_{pq} m_{ra}=
\sum_i \Bigl( \sum_a x_{jl}^{ak}x_{aq}^{ip} \Bigl ) m_{ri}$$
and
$$(m_{pq}\cdot m_{kl})\cdot m_{rj}=\sum_a x_{lq}^{ap}m_{ka}m_{rj}=
\sum_i \Bigl(\sum_a x_{ja}^{ik}x_{lq}^{ap} \Bigl)m_{ri}$$
such that
$$
(m_{pq}\cdot m_{kl})\cdot m_{rj}-m_{pq}\cdot (m_{kl}\cdot m_{rj})=
\sum_i \Bigl(\sum_a x_{ja}^{ik}x_{lq}^{ap} - x_{jl}^{ak}x_{aq}^{ip} \Bigl)
m_{ri}
$$
and the right hand side is zero for all indices $p,q,k,l,r$ and $j$
if and only if \eqref{3.3.1} holds, and in \prref{3.3}, we have
seen that \eqref{3.3.1} is equivalent to $R^{12}R^{23}=R^{13}R^{12}$.\\
The last statement
follows from
$$(m_{pq}\cdot m_{kl})\bullet m_{j}-m_{pq}\bullet (m_{kl}\bullet m_{j})=
\sum_i \Bigl(\sum_a x_{ja}^{ik}x_{lq}^{ap} - x_{jl}^{ak}x_{aq}^{ip} \Bigl)
m_{i}=0$$
where we used \eqref{3.3.1} at the last step.
\end{proof}

\begin{examples}\rm\exlabel{3.5}
1) The identity $I_{M\ot M}$ and the switch map $\tau_M$ are
trivial solutions of the FS-equation in $\End_k(M\ot M\ot M)$.\\

2) Let $A={\cal M}_n(k)$ and $(e_{i}^{j})_{1\leq i,j\leq n}$ be the usual
basis. Then
$$R=\sum_{i=1}^n e_{i}^{j}\ot e_{j}^{i}$$
is  a solution of the $S$-equation and
$(R,{\rm trace})$ is not a solution of the F-equation.\\

3) Let $R\in A\ot A$ be a solution of FS-equation and
$u\in A$ invertible. Then
\begin{equation}\eqlabel{3.5.1}
{}^uR=(u\ot u)R(u^{-1}\ot u^{-1})
\end{equation}
is also a solution
of the FS-equation. Let ${\bf FS}(A)$ be the set of all
solutions of the FS-equation, and $U(A)$ the multiplicative group
of invertible elements in $A$. Then \eqref{3.5.1} defines an
action of $U(A)$ on ${\bf FS}(A)$.\\

4) If $a\in A$ is an idempotent, then $a\ot a$ is a solution of the
FS-equation.\\

5) Let $A$ be a $k$-algebra, and $e\in A\ot A$ an $A$-central element.
Then for any left $A$-module $M$, the map
$R=R_e:\ M\ot M\to M\ot M$ given by
\begin{equation}\eqlabel{2.4.1}
R(m\ot n)=\sum e^1\cdot m\ot e^2\cdot n
\end{equation}
$(m,n\in M)$ is a solution of the FS-equation in $\End_k(M\ot M\ot M)$.
This is an easy consequence of \eqref{2.3.2a}.
Moreover, if $e$ is a separability idempotent
(respectively $(e,\eps)$ is a Frobenius pair),
then $R$ is a solution of the S-equation (respectively a solution of the
F-equation).\\

6) Let $G$ be a finite group, and $A=kG$. Then
$e=\displaystyle{\sum_{g\in G}} g\ot g^{-1}$ is an $A$-central element and
$(e,p_1)$ is a Frobenius pair ($p_1:\ kG\to k$ is the map defined by
$p_1(g)=\delta_{1,g}$, for all $g\in G$). Hence,
$kG$ is a Frobenius algebra (this is the Frobenius reciprocity law from
group representation theory).
Furthermore, if $({\rm char}(k), |G|)=1$, then
$e^{\prime}=|G|^{-1}e$ is a separability idempotent.\\

7) Using a computer program, Bogdan Ichim computed for us that
${\bf FS}({\cal M}_2(\nint_2))$ (resp. ${\bf FS}({\cal M}_2(\nint_3))$)
consists of exactly 38 (resp. 187) solutions of the FS-equation.
We will present only two of them.
Let $k$ be a field of characteristic $2$ (resp. $3$). Then
$$
R=\pmatrix{
1&1&1&0\cr
1&0&0&1\cr
1&0&0&1\cr
0&1&1&1\cr}, \quad ({\rm resp.} \quad
R=\pmatrix{
1&0&0&1\cr
0&1&1&2\cr
0&1&1&2\cr
1&2&2&1\cr} \;\; )
$$
are solutions of $FS$-equation.\\

8) Let $M$ be a free $k$-module, with basis $\{m_1,m_2,\cdots,m_n\}$,
and let $R$ be given by
$$R(m_u\ot m_v)=\sum_j a^j_u m_v\ot m_j$$
Thus $x^{ij}_{uv}=\delta^i_v a^j_u$, with notation as in \eqref{matrix5}.
An immediate verification shows that $R$ is a solution of the
FS-equation. If $\varepsilon$ is the trace map, then $(R,\varepsilon)$
is a solution of the F-equation if and only if $a^j_i=\delta^j_i$.
$R$ is a solution of the $S$-equation if and only if $n$ is invertible in
$k$, and $na^j_i=\delta^j_i$.
\end{examples}

If $H$ is a finite dimensional unimodular involutory Hopf algebra,
and $t$ is a two-sided integral in $H$, then
$R=\sum t_{(1)}\ot S(t_{(2)})$ is a solution of the
quantum Yang-Baxter equation (cf. \cite[Theorem 8.3.3]{LR}).
In our next Proposition, we will show that, for an arbitrary Hopf algebra $H$,
$R$ is a solution of the
FS-equation and the braid
equation.

\begin{proposition}\prlabel{3.6}
Let $H$ be a Hopf algebra over a commutative ring $k$, and $t\in H$
a left integral. Then $R=\sum t_{(1)}\ot S(t_{(2)})\in H\ot H$
is $H$-central, and therefore
a solution of the FS-equation and the braid equation.
\end{proposition}

\begin{proof}
For all $h\in H$, we have that $ht=\varepsilon(h)t$ and, subsequently,
$$\sum h_{(1)}t\ot h_{(2)}=t\ot h$$
$$\sum h_{(1)}t_{(1)}\ot S(h_{(2)}t_{(2)})\ot h_{(3)}=
\sum t_{(1)}\ot S(t_{(2)})\ot h$$
Multiplying the second and the third factor, we obtain
$$\sum ht_{(1)}\ot S(t_{(2)})=\sum t_{(1)}\ot S(t_{(2)})h$$
proving that $R$ is $H$-central.
\end{proof}

\begin{remarks}\relabel{3.6.b}\rm
1) If $\varepsilon(t)=1$, then $R$ is a separability idempotent and
$H$ is separable over $k$, and, using \reref{1.1}, we recover
Maschke's Theorem for Hopf algebras (see \cite{Ab}).\\
2) If $t$ is a right integral, then a similar argument shows that
$\sum S(t_{(1)})\ot t_{(2)}$ is an $H$-central element.
\end{remarks}

Let $H$ be a Hopf algebra and $\sigma:\ H\ot H\to k$ be a normalized
Sweedler 2-cocycle, i.e.
$$\sigma (h\ot 1_H)=\sigma (1_H\ot h)=\eps (h)$$
\begin{equation}\eqlabel{3.7.1c}
\sum \sigma (k_{(1)}\ot l_{(1)})
\sigma (h\ot k_{(2)}l_{(2)})= \sum \sigma (h_{(1)}\ot k_{(1)})
\sigma (h_{(2)}k_{(2)}\ot l)\end{equation}
for all $h$, $k$, $l\in H$. The
{\sl crossed product} algebra $H_{\sigma}$ is equal to $H$
as $k$-module and the (associative) multiplication is given by
$$g\cdot h = \sum \sigma (g_{(1)}\ot h_{(1)})g_{(2)}h_{(2)}$$
for all $g$, $h\in H_{\sigma}=H$.

\begin{proposition}\prlabel{3.7}
Let $H$ be a cocommutative Hopf algebra over a commutative ring $k$,
$t\in H$ a right integral,
and $\sigma:\ H\ot H\to k$ a normalized
convolution invertible Sweedler 2-cocycle. Then
\begin{equation}\eqlabel{tsigm}
R_{\sigma}= \sum \sigma^{-1} \Bigl( S(t_{(2)})\ot t_{(3)} \Bigl)
S(t_{(1)})\ot t_{(4)}\in H_{\sigma}\ot H_{\sigma}
\end{equation}
is $H_{\sigma}$-central, and
a solution of the FS-equation and the braid equation.
Consequently, if $H$ is a separable Hopf algebra, then $H_{\sigma}$
is also a separable algebra.
\end{proposition}

\begin{proof}
The method of proof is the same as in \prref{3.6}, but the situation
is more complicated. For all $h\in H$, we have $th=\varepsilon(h)t$,
and
$$h\ot t=\sum h_{(1)}\ot th_{(2)}$$
and
$$\sum h\ot S(t_{(1)})\ot S(t_{(2)})\ot t_{(3)}\ot t_{(4)}=
\sum h_{(1)}\ot S(h_{(2)})S(t_{(1)})\ot S(h_{(3)})S(t_{(2)})
\ot t_{(3)}h_{(4)}\ot t_{(4)}h_{(5)}$$
We now compute easily that
\begin{eqnarray*}
h\cdot R_{\sigma}&=&
\sum \sigma^{-1}(S(t_{(2)})\ot t_{(3)})h\cdot S(t_{(1)})\ot t_{(4)}\\
&=&\sum \sigma^{-1}(S(h_{(3)})S(t_{(2)}\ot t_{(3)}h_{(4)})
h_{(1)}\cdot(S(h_{(2)})S(t_{(1)}))\ot t_{(4)}h_{(5)}\\
&=& \sum \sigma^{-1}(S(h_{(3)})S(t_{(2)}\ot t_{(3)}h_{(4)})
\sigma((h_{(1)})_{(1)}\\
&&~~~~~~~~\ot S(h_{(2)})_{(1)}S(t_{(1)})_{(1)})
(h_{(1)})_{(2)}S(h_{(2)})_{(2)}S(t_{(1)})_{(2)}\ot t_{(4)}h_{(5)}\\
&=& \sum \sigma^{-1}(S(h_{(3)})S(t_{(3)}\ot t_{(4)}h_{(4)})
\sigma(h_{(1)}\ot S(h_{(2)}S(t_{(2)})) S(t_{(1)})\ot t_{(5)}h_{(5)}
\end{eqnarray*}
On the other hand
\begin{eqnarray*}
R_{\sigma}\cdot h&=&
\sum \sigma^{-1}(S(t_{(2)})\ot t_{(3)})S(t_{(1)})\ot t_{(4)}\cdot h\\
&=& \sum \sigma^{-1}(S(t_{(2)})\ot t_{(3)})\sigma(t_{(4)}\ot h_{(1)})
S(t_{(1)})\ot t_{(5)}h_{(2)}
\end{eqnarray*}
In order to prove that $R_{\sigma}$ is $H_{\sigma}$-central, it suffices
to show that, for all $f,g\in H$:
\begin{equation}\eqlabel{3.7.1}
\sum \sigma^{-1}(S(h_{(3)})S(g_{(2)})\ot g_{(3)}h_{(4)})
\sigma(h_{(1)}\ot S(h_{(2)})S(g_{(1)}))=
\sum \sigma^{-1}(S(g_{(1)})\ot g_{(2)})\sigma(g_{(3)}\ot h)
\end{equation}
So far, we have not used the cocommutativity of $H$. If $H$ is cocommutative,
then we can omit the Sweedler indices, since they contain no information.
Hence we can write
$$\Delta(h)=\sum h\ot h$$
The cocycle relation \eqref{3.7.1c} can then be rewritten as
\begin{eqnarray}
\sum \sigma(h\ot kl)&=& \sum \sigma^{-1}(k\ot l)\sigma(h\ot k)
\sigma(kh\ot l)\eqlabel{3.7.1a}\\
\sum \sigma(hk\ot l)&=& \sum \sigma^{-1}(h\ot k)\sigma(k\ot l)
\sigma(h\ot kl)\eqlabel{3.7.1b}
\end{eqnarray}
Using \eqref{3.7.1a}, \eqref{3.7.1b} and the fact that $\sigma$ is
normalized, we compute
\begin{eqnarray*}
&&\hspace*{-2cm} \sum \sigma^{-1}(S(h)S(g)\ot gh)
\sigma(h\ot S(h)S(g))\\
&=& \sum \sigma(S(h)\ot S(g))\sigma^{-1}(S(g)\ot gh)
\sigma^{-1}(S(h)\ot S(g)gh)\\
&&\hspace*{1cm}\sigma^{-1}(S(h)\ot S(g))
\sigma(h\ot S(h))\sigma(hS(h)\ot g)\\
&=& \sum \sigma(g\ot h)\sigma^{-1}(S(g)\ot g)\sigma^{-1}(S(g)g\ot h)\\
&=& \sum \sigma^{-1}(S(g)\ot g)\sigma(g\ot h)
\end{eqnarray*}
proving \eqref{3.7.1}. We also used that
$$\sum \sigma(h\ot S(h))=\sum \sigma(S(h)\ot h)$$
which follows from the cocycle condition and the fact that $\sigma$
is normalized. Finally, if $H$ is separable, then we can find a right
integral $t$ such that $\varepsilon(t)=1$, and we easily see that
$m_{H_{\sigma}}(R_{\sigma})=1$, proving that $R_{\sigma}$ is a solution
of the $S$-equation.
\end{proof}

In \cite{ESS} solutions of the braid equation are constructed starting
from $1$-cocycles on a group $G$. The interesting point in this
construction is that, at set theory level, any ``nondegenerate
symmetric" solution of the braid equation arises in this way
(see \cite[Theorem 2.9]{ESS}). Now, taking $G$ a finite group and
$H=kG$ in \prref{3.7}, we find a large class of solutions to
the braid equation, arising from $2$-cocycles $\sigma:\ G\times G\to k^*$.
These solutions $R$ can be described
using a family of
scalars $(x_{uv}^{ij})$, as in \prref{3.3}, where the indices now run through
$G$. Let $n=|G|$, and write
${\cal M}_n(k)\cong {\cal M}_G(k)$, with entries indexed by $G\times G$.

\begin{corollary}\colabel{3.9}
Let $G$ be a finite group of order $n$, and $\sigma:\ G\times G\to k^*$ a
normalized $2$-cocycle. Then $R_{\sigma}=(x_{uv}^{ij})_{i,j,u,v\in G}$ given by
\begin{equation}\eqlabel{cocsol}\eqlabel{3.9.1}
x_{uv}^{ij}=\delta_{j,\; ui^{-1}v}\;\; \sigma^{-1}(iu^{-1}, ui^{-1})
\sigma (iu^{-1}, u) \sigma (ui^{-1}, v)
\end{equation}
($i$, $j$, $u$, $v\in G$) is a solution of the FS-equation.
If $n$ is invertible in $k$, then $n^{-1}R$ is a solution of the
S-equation.
\end{corollary}

\begin{proof}
The twisted group algebra $k_{\sigma}G$ is the $k$-module with basis $G$,
and multiplication given by $g\cdot h=\sigma (g,h)gh$, for any
$g$, $h\in G$. $t=\sum_{g\in G}g$ is a left integral in $kG$ and the
element $R_{\sigma}$ defined in \eqref{tsigm} takes the form
$$R_{\sigma}= \sum_{g\in G}\; \sigma^{-1} (g^{-1}, g)g^{-1}\ot g$$
Using the multiplication rule on $k_{\sigma}G$, we find that the map
$$\tilde{R_{\sigma}}: k_{\sigma}G\ot k_{\sigma}G\to k_{\sigma}G\ot
k_{\sigma}G, \quad
\tilde{R_{\sigma}} (u\ot v)= R_{\sigma}\cdot (u\ot v)$$
is given by
$$\tilde{R_{\sigma}} (u\ot v)=\sum_{i\in G}\;
\sigma^{-1}(iu^{-1}, ui^{-1})\sigma (iu^{-1}, u)\sigma (ui^{-1}, v)\;
i\ot ui^{-1}v$$
If we write
$$\tilde{R_{\sigma}} (u\ot v)=\sum_{i,j\in G}\; x_{uv}^{ij}\; i\ot j$$
then $x_{uv}^{ij}$ is given by \eqref{cocsol}.
\end{proof}

We will now present a coalgebra version of \exref{3.5}, 3).
First we adapt an old definition of Larson (\cite{L}).

\begin{definition}\delabel{3.10}
Let $C$ be a $k$-coalgebra.
A map $\sigma:\ C\ot C\to k$ is called an FS-map if
\begin{equation}\eqlabel{3.10.1}
\sum \sigma(c\ot d_{(1)})d_{(2)}=\sum \sigma(c_{(2)}\ot d)c_{(1)}.
\end{equation}
If, in addition, $\sigma$ satisfies the normalizing condition
\begin{equation}\eqlabel{3.10.2}
\sum \sigma(c_{(1)}\ot c_{(2)})=\varepsilon(c)
\end{equation}
then $\sigma$ is called a coseparability idempotent.\\
If there exists an $f\in C$ such that the FS-map $\sigma$ satisfies
the normalizing condition
\begin{equation}\eqlabel{3.10.3}
\sigma(f\ot c)=\sigma(c\ot f)=\varepsilon(c)
\end{equation}
for all $c\in C$, then we call $(\sigma,f)$ an F-map.
\end{definition}

\begin{proposition}\prlabel{3.10b}
Let $C$ be a coalgebra over a commutative ring $k$.\\
1) $C$ is coseparable if and only if there exists a coseparability
idempotent $\sigma$;\\
2) The corestriction of scalars
(or forgetful) functor ${\cal M}^C\to {}_k{\cal M}$
is Frobenius
 if and only if there exists an $F$-map $(\sigma,f)$.
\end{proposition}

\begin{proof}
1) is well-known, see \cite{L}. To prove 2), we recall first that
$G=C\ot \cdot$ is a right adjoint of the forgetful functor
$F:\ {\cal C}={}^C{\cal M}\to {\cal M}={}_k{\cal M}$. $G$ is also a left
adjoint
of $F$ if there exists natural transformations
$$\theta:\ 1_{\cal M}\to FG~~{\rm and}~~\nu:\ GF\to 1_{\cal C}$$
such that
\begin{eqnarray}
F(\nu_M)\circ \theta_{F(M)}&=& I_M\eqlabel{3.10b.1}\\
\nu_{G(N)}\circ G(\theta_N)&=& I_N\eqlabel{3.10b.2}
\end{eqnarray}
for all $M\in {}^C{\cal M}$, $N\in {\cal M}$. In \cite{CIMZ}, we have
seen that there is a one-to-one correspondence between natural
transformations $\nu:\ GF\to 1_{\cal C}$ and FS-maps
(see \cite[Theorem 2.3, Proposition 2.5 and Example 2.11]{CIMZ}).
Given $\nu$, the corresponding FS-map is $\sigma=\varepsilon_C\circ\nu_C$.
Given $\sigma$, we recover $\nu$ by putting $\nu_M:\ GF(M)=C\ot M\to M$,
$$\nu_M(c\ot m)=\sum \sigma(c\ot m_{(-1)})m_{(0)}$$
There is also a one-to-one correspondence between natural transformations
$\theta:\ 1_{\cal M}\to FG$, and elements $f\in C$
(see \cite[Theorem 2.4]{CMZ}, in the special case $A=k$). The element
$f\in C$ corresponding to $\theta$ is just $\theta_k(1)$. Conversely,
given $f$, we construct $\theta$ by setting $\theta_N:\ N\to C\ot N$,
$$\theta_N(n)=c\ot n$$
for every $k$-module $N$. \eqref{3.10b.1} is equivalent to
\begin{equation}\eqlabel{3.10b.3}
\sum \sigma(f\ot m_{(-1)})m_{(0)}=m
\end{equation}
for all $m\in M$. Taking $M=C$ and $m=c$, we obtain, after applying
$\varepsilon$ to both sides
\begin{equation}\eqlabel{3.10b.4}
\sigma(f\ot c)=\varepsilon(c)
\end{equation}
If \eqref{3.10b.4} holds, then we find for all $m\in M$
$$\sum \sigma(f\ot m_{(-1)})m_{(0)}=\sum \varepsilon(m_{(-1)})m_{(0)}=m$$
and \eqref{3.10b.3} holds.\\
Now \eqref{3.10b.2} is equivalent to
\begin{equation}\eqlabel{3.10b.5}
\sum \sigma(c\ot f_{(1)})f_{(2)}\ot n=c\ot n
\end{equation}
for all $c\in C$ and $n\in N$. Taking $N=k$, $n=1$, and applying
$\varepsilon$ to both sides, we find
\begin{equation}\eqlabel{3.10b.6}
\sigma(c\ot f)=\varepsilon(c)
\end{equation}
for all $c\in C$. Conversely, if \eqref{3.10b.6} holds, then
\begin{eqnarray*}
\sum \sigma(c\ot f_{(1)})f_{(2)}\ot n&=&
\sum \sigma(c_{(2)}\ot f)c_{(1)}\ot n\\
\sum \varepsilon(c_{(2)})c_{(1)}\ot n&=&c\ot n
\end{eqnarray*}
and \eqref{3.10b.5} follows.
\end{proof}

\begin{examples}\relabel{3.11}\rm
1) The comatrix coalgebra ${\cal M}^n(k)$ is coseparable and
$$\sigma:\ {\cal M}^n(k)\ot {\cal M}^n(k) \to k, \quad
\sigma (c^i_{j}\ot c^k_{l})=\delta_{kj}\delta_{il}$$
is a coseparability idempotent.\\

2) Let $k$ be a field of characteristic zero, and consider the
Hopf algebra
$C=k[Y]$ (with $\Delta (Y)=Y\ot 1 +1\ot Y$ and $\varepsilon(Y)=0$).
Then the only FS-map $\sigma:\ C\ot C\to k$ is the zero map.
Indeed, $sigma$ is completely described by $\sigma(Y^i\ot Y^j)=a_{ij}$.
We will show that all $a_{ij}=0$. Using the fact that $\Delta(Y^n)=
\delta(Y)^n$, we easily find that \eqref{3.10.1} is equivalent to
\begin{equation}\eqlabel{3.112).1}
\sum_{j=0}^m {m \choose j}a_{n,m-j}Y^j=
\sum_{i=0}^n {n \choose i}a_{n-i,m}Y^i
\end{equation}
for all $n,m\in \nnat$. Taking $n>m$, and identifying the coefficients
in $Y^n$, we find $a_{n,m}=0$. If $m>n$, we also find $a_{n,m}=0$,
now identifying coefficients in $Y^m$. We can now write
$a_{nm}=a_n\delta_{nm}$. Take $m>n$. The righthand side of
\eqref{3.112).1} amounts to zero, while the lefthand side is
$${m \choose m-n}a_nY^{n-m}$$
It follows that $a_n=0$ for all $n$, and $\sigma=0$.
\end{examples}

\begin{proposition}\prlabel{3.12}
Let $C$ be a coalgebra, $\sigma:\ C\ot C\to k$ an $FS$-map
and $M$ a right $C$-comodule. Then the map
$$R_{\sigma}:\ M\ot M\to M\ot M, \quad
R_{\sigma}(m\ot n)=\sum \sigma(m_{(1)}\ot n_{(1)})m_{(0)}\ot n_{(0)}$$
($m$, $n\in M$) is a solution of the $FS$-equation in
$\End_k(M\ot M\ot M)$.
\end{proposition}

\begin{proof}
Write $R=R_{\sigma}$ and take $l$, $m$, $n\in M$.
\begin{eqnarray}
R^{12}R^{23}(l\ot m\ot n)&=&
R^{12}\Bigl(\sum \sigma(m_{(1)}\ot n_{(1)})
l\ot m_{(0)}\ot n_{(0)}\Bigl)\nonumber\\
&=&\sum \sigma(m_{(2)}\ot n_{(1)})
\sigma(l_{(1)}\ot m_{(1)})l_{(0)}\ot m_{(0)}\ot n_{(0)}.\eqlabel{2.6.1}
\end{eqnarray}
Applying \eqref{3.10.1} to \eqref{2.6.1}, with $m=c,~n=d$, we obtain
\begin{eqnarray*}
R^{12}R^{23}(l\ot m\ot n)&=&
\sum \sigma(m_{(1)}\ot n_{(1)})\sigma(l_{(1)}\ot n_{(2)})
l_{(0)}\ot m_{(0)}\ot n_{(0)}\\
&=& R^{23}\Bigl(\sum \sigma(l_{(1)}\ot n_{(1)})
l_{(0)}\ot m\ot n_{(0)}\\
&=& R^{23}R^{13}(l\ot m\ot n).
\end{eqnarray*}
Applying \eqref{3.10.1} to \eqref{2.6.1}, with $m=d,~l=c$, we obtain
\begin{eqnarray*}
R^{12}R^{23}(l\ot m\ot n)&=&
\sum \sigma(l_{(1)}\ot n_{(1)})\sigma(l_{(2)}\ot m_{(1)})
l_{(0)}\ot m_{(0)}\ot n_{(0)}\\
&=& R^{13}\Bigl(\sum \sigma(l_{(1)}\ot n_{(1)})
l_{(0)}\ot m_{(0)}\ot n\\
&=& R^{13}R^{12}(l\ot m\ot n)
\end{eqnarray*}
proving that $R$ is a solution of the $FS$-equation.
\end{proof}

\begin{remark}\relabel{3.13}\rm
If $C$ is finitely generated and projective, and $A=C^*$ is its dual
algebra, then
there is a one-to-one correspondence between $FS$-maps
$\sigma:\ C\ot C\to k$ and $A$-central elements $e\in A\ot A$.
The correspondence is given by the formula
$$\sigma(c\ot d)=\sum \lan c,e^1\ran \lan d,e^2\ran .$$
In this situation, the map $R_e$ from \exref{3.5} 5) is equal to
$R_{\sigma}$. Indeed,
\begin{eqnarray*}
R_{\sigma}(m\ot n)&=&
\sum \sigma(m_{(1)}\ot n_{(1)}) m_{(0)}\ot n_{(0)}\\
&=&\sum \lan m_{(1)},e^1\ran \lan n_{(1)},e^2\ran\\
&=& \sum e^1\cdot m\ot e^2\cdot n\\
&=& R_e(m\ot n).
\end{eqnarray*}
\end{remark}

We will now present two more
classes of solutions of the $FS$-equation.

\begin{proposition}\prlabel{lazio}\prlabel{3.14}
Take $a\in k$, $X=\{1,\dots,n\}$, and $\theta:\ X^3\to X$
a map satisfying
\begin{eqnarray}
\theta(u,i,j)=v&\Longleftrightarrow& \theta(v,j,i)=u\eqlabel{3.14.1}\\
\theta(i,u,j)=v&\Longleftrightarrow& \theta(j,v,i)=u\eqlabel{3.14.2}\\
\theta(i,j,u)=v&\Longleftrightarrow& \theta(j,i,v)=u\eqlabel{3.14.3}\\
\theta(i,j,k)=\theta(u,v,w)&\Longleftrightarrow& \theta(j,i,u)=
\theta(k,w,v)\eqlabel{3.14.4}
\end{eqnarray}
1) $R=(x^{uv}_{ij})\in {\cal M}_{n^2}(k)$ given by
\begin{equation}\eqlabel{3.14.4bis}
x^{uv}_{ij}=a\delta^u_{\theta(i,v,j)}
\end{equation}
is a solution of the FS-equation.\\
2) Assume that $n\in \nnat$ is invertible in $k$, and take $a=n^{-1}$.\\
2a) $R$ is a solution of the S-equation if and only if
$$\theta(k,k,i)=i$$
for all $i,k\in X$.\\
2b) Let $\varepsilon$ be the trace map. $(R,\varepsilon)$ is a solution of
the F-equation if and only if
$$\theta(i,k,k)=i$$
for all $i,k\in X$.
\end{proposition}

\begin{proof}
1) We have to verify \eqref{3.3.1}. Using \eqref{3.14.3}, we compute
\begin{eqnarray}
\sum_k x^{ij}_{uk}x^{kl}_{vw}&=&
a^2\sum_k \delta^i_{\theta(u,j,k)}\delta^k_{\theta(v,l,w)}\nonumber\\
&=&a^2\sum_k \delta^k_{\theta(j,u,i)}\delta^k_{\theta(v,l,w)}\nonumber\\
&=&a^2\sum_k \delta^{\theta(j,u,i)}_{\theta(v,l,w)}\eqlabel{3.14.5}
\end{eqnarray}
In a similar way, we find
\begin{eqnarray}
\sum_k x^{jl}_{vk}x^{ik}_{uw}&=&
a^2\sum_k \delta^{\theta(w,i,u)}_{\theta(l,v,j)}\eqlabel{3.14.6}\\
\sum_k x^{il}_{kw}x^{kj}_{uv}&=&
a^2\sum_k \delta^{\theta(i,w,l)}_{\theta(u,j,v)}\eqlabel{3.14.7}
\end{eqnarray}
Using \eqref{3.14.4}, we find that \eqref{3.14.5}, \eqref{3.14.6},
and \eqref{3.14.7} are equal, proving \eqref{3.3.1}.\\

2a) We easily compute that
$$\sum_k x^{kj}_{ik}=n^{-1}\sum_k \delta^k_{\theta(i,j,k)}=
n^{-1}\sum_k \delta^j_{\theta(k,k,i)}$$
and it follows from  that $R$ is a solution of the S-equation if and
only if $\theta(k,k,i)=i$ for all $i$ and $k$.\\

2b) We compute
\begin{eqnarray*}
\sum_k x^{kj}_{ki}&=& n^{-1}\sum_k \delta^k_{\theta(k,j,i)}=
n^{-1}\sum_k \delta^j_{\theta(i,k,k)}\\
\sum_k x^{jk}_{ik}&=&n^{-1}\sum_k \delta^j_{\theta(i,k,k)}
\end{eqnarray*}
and it follows from  that
$(R,\varepsilon)$ is a solution of
the F-equation if and only if
$\theta(i,k,k)=i$
for all $i,k$.
\end{proof}

\begin{examples}\exlabel{3.15}\rm
1) Let $G$ be a finite group. Then the map
$$\theta:\ G\times G\times G\to G,~~\theta(i,j,k)=ij^{-1}k$$
satisfies conditions (\ref{eq:3.14.1}-\ref{eq:3.14.4}).\\

2) Let $G$ be a group of order $n$ acting on $X=\{1,2,\cdots, n\}$,
and assume that the action of $G$ is transitive and free, which
means that for every $i,j\in X$, there exists a unique $g\in G$
such that $g(i)=j$. Then the map
$\theta:\ X\times X\times X\to X$ defined by
$$\theta(i,j,k)=g^{-1}(k)$$
where $g\in G$ is such that $g(i)=j$,
satisfies conditions (\ref{eq:3.14.1}-\ref{eq:3.14.4}).
\end{examples}

\begin{proposition}\prlabel{3.16}
Let $n$ be a positive integer, $\phi: \{1, \cdots, n \} \to
\{1, \cdots, n \}$ a function with $\phi^{2}=\phi$
and $M$ a free of rank $n$ $k$-module with
$\{m_1,\cdots, m_n\}$ a $k$-basis. Then
\begin{equation}\eqlabel{3.16.1}
R^{\phi}:\ M\ot M\to M\ot M,\quad
R^{\phi}(m_i\ot m_j)=\delta_{ij}
\sum_{a,b\in \phi^{-1}(i)}m_a\ot m_b
\end{equation}
($i$, $j=1,\cdots, n$), is a solution of the $FS$-equation.
\end{proposition}

\begin{proof} Write $R=R^{\phi}$, and take $p,q,r\in
\{1,\ldots,n\}$.
Then
\begin{eqnarray*}
R^{12}R^{23}(m_p\ot m_q\ot m_r)&=&
R^{12}\Bigl (\delta_{qr}
\sum_{a,b\in \phi^{-1}(q)}m_p\ot m_a\ot m_b \Bigl )\\
&=&\delta_{qr}\sum_{a,b\in \phi^{-1}(q)}
\delta_{ap}\sum_{c,d\in \phi^{-1}(p)}
m_c\ot m_d\ot m_b\\
&=&\delta_{qr}\delta_{\phi(p)q}\sum_{b\in \phi^{-1}(q)}
\sum_{c,d\in \phi^{-1}(p)} m_c\ot m_d\ot m_b.
\end{eqnarray*}
If $\phi^{-1}(p)$ is nonempty (take $x\in \phi^{-1}(p)$), and $\phi(p)=q$,
then $q=\phi(p)=\phi^2(x)
=\phi(x)=p$, so we can write
$$R^{12}R^{23}(m_p\ot m_q\ot m_r)=
\delta_{pqr\phi(p)} \sum_{a,b,c\in \phi^{-1}(p)} m_a\ot m_b\ot m_c.$$
In a similar way, we can compute that
\begin{eqnarray*}
R^{23}R^{13}(m_p\ot m_q\ot m_r)&=&
R^{13}R^{12}(m_p\ot m_q\ot m_r)\\
&=&\delta_{pqr\phi(p)} \sum_{a,b,c\in \phi^{-1}(p)} m_a\ot m_b\ot m_c.
\end{eqnarray*}
\end{proof}

Note that $R^{\phi}$ appeared already in \cite{M4}, where it
was proved that $R^{\phi}$ is a solution of the integrability condition of
the Knizhnik-Zamolodchikov equation $[R^{12}, R^{13}+R^{23}]=0$.

Now we will generalize
\exref{3.5} 3) to algebras
without a unit. Recall that a left $A$-module $M$ is called unital
(or regular, in the terminology of \cite{T}) if the natural map
$A\ot_A M\to M$ is an isomorphism.

\begin{proposition}\prlabel{3.15}
Let $M$ be a unital $A$-module, and $f:\ A\to A\ot A$ an $A$-bimodule map.
Then the map $R:\ M\ot M\to M\ot M$ mapping $m\ot a\cdot n$ to
$f(a)(m\ot n)$ is a solution of the $FS$-equation.
\end{proposition}

\begin{proof}
Observe first that it suffices to define $R$ on elements of the form
$m\ot a\cdot n$, since the map $A\ot_A M\to M$ is surjective.
$R$ is well-defined since
\begin{eqnarray*}
R(m\ot a\cdot (b\cdot n))&=& f(a)(m\ot b\cdot n)\\
&=& f(a)(I_M\ot b)(m\ot n)\\
&=& f(ab)(m\ot n)
\end{eqnarray*}
Write $f(a)=\sum a^1\ot a^2$, for all $a\in A$. Then
$$f(ab)=\sum a^1\ot a^2b=\sum ab^1\ot b^2.$$
Now
\begin{eqnarray*}
R^{12}(R^{23}(m\ot bn \ot ap))&=&
R^{12}(\sum m\ot a^1bn\ot a^2 p)\\
&=& \sum a^1b^1m\ot b^2n\ot a^2p\\
&=& R^{13} (\sum b^1m\ot b^2n\ot p)\\
&=& R^{13}(R^{12}(m\ot bn \ot ap)).
\end{eqnarray*}
In a similar way, we prove that $R^{12}R^{23}=R^{23}R^{13}$.
\end{proof}

\section{The structure of Frobenius algebras and
separable algebras}\selabel{4}
The first statement of \prref{nori} can be restated as follows:
for a $k$-algebra $A$, any $A$-central element $R\in A\ot A$
is a solution of the FS-equation.
We will now prove the converse: for a flat $k$-algebra $A$,
any solution of the FS-equation arises in
this way.

\begin{theorem}\thlabel{4.1}
Let $A$ be a flat $k$-algebra and
$R=\sum R^1\ot R^2\in A\ot A$
a solution of the FS-equation. Then
\begin{enumerate}
\item[1)] Then there exists a $k$-subalgebra
${\cal A}(R)$ of $A$ such that $R\in {\cal A}(R)\ot {\cal A}(R)$
and $R$ is ${\cal A}(R)$-central.
\item[2)] If $R\sim S$ are equivalent solutions of the $FS$-equation,
then ${\cal A}(R)\cong {\cal A}(S)$.
\item[3)] $({\cal A}(R), R)$ satisfies the following universal property:
if $(B, m_B, 1_B)$ is an algebra, and $e\in B\ot B$
is an $B$-central element, then any algebra map
$\alpha:\ B\to A$ with $(\alpha\ot\alpha)(e)=R$ factors through an algebra
map $\tilde{\alpha}:\ B\to {\cal A}(R)$.
\item[4)] If $R\in A\ot A$ is a solution of the S-equation
(resp. the F-equation), then $\cal{A}(R)$ is a separable
(resp. Frobenius) algebra.
\end{enumerate}
\end{theorem}

\begin{proof} 1) Let
${\cal A}(R)=\{a\in A~|~a\cdot R=R\cdot a\}$.
Obviously ${\cal A}(R)$ is a $k$-subalgebra of $A$ and
$1_A\in {\cal A}(R)$.
We also claim that $R\in {\cal A}(R)\ot {\cal A}(R)$.
To this end, we observe first that
${\cal A}(R)={\rm Ker}(\varphi)$, with
$\varphi:\ A\to A\ot A^{\rm op}$ defined by the formula
$$\varphi(a)=(a\ot 1_A-1_A\ot a)R.$$
$A$ is flat as a $k$-algebra, so
$${\cal A}(R)\ot A= {\rm Ker}(\varphi\ot {\rm Id}_{A}).$$
Now,
\begin{eqnarray*}
(\varphi\ot I_{A})(R)&=&
\sum r^1R^1\ot R^2\ot r^2- \sum R^1\ot R^2r^1\ot r^2\\
&=& R^{13}R^{12}-R^{12}R^{23}=0
\end{eqnarray*}
and it follows that $R\in {\cal A}(R)\ot A$. In a similar way,
using that $R^{12}R^{23}=R^{23}R^{13}$, we get that
$R\in A\ot {\cal A}(R)$, and we find that
$R\in {\cal A}(R)\ot {\cal A}(R)$. Hence, $R$ is an ${\cal A}(R)$-central
element of ${\cal A}(R)\ot {\cal A}(R)$.\\
2) Let $u\in U(A)$ such that $S=(u\ot u)R(u^{-1}\ot u^{-1})$. Then
$$f_u :{\cal A}(R) \to {\cal A}(S), \quad
f_u (a)=uau^{-1}$$
for all $a\in {\cal A}(R)$ is a well-defined isomorphism of $k$-algebras.\\
3) Let $b\in B$. If we apply $\alpha\ot \alpha$ to the equality
$(b\ot 1_B)e=e(1_B\ot b)$ we find that the image of $\alpha$ is
contained in ${\cal A}(R)$, and the universal property follows.\\
4) The first statement follows from the definition of
separable algebras and the second one follows from 3) of
\coref{adi}.
\end{proof}

\begin{remark}
\rm Let $A$ be a separable $k$-algebra and $R\in A\ot A$ a
separability idempotent. Then $R$ and $S:=0\ot 0$ are solutions of the
$FS$-equation, ${\cal A}(R)={\cal A}(S)=A$ and of course, $R$ and $S$ are
not equivalent.
\end{remark}

We can now prove the main result of this paper.

\begin{theorem}\thlabel{imp}\thlabel{4.2}
For a finitely generated and projective algebra $A$
over a commutative ring $k$, the following statements are
equivalent:
\begin{enumerate}
\item[1)] $A$ is a Frobenius (resp. separable) algebra.
\item[2)] There exists an algebra isomorphism
$$A\cong {\cal A}(R),$$
where $R\in \End_k(A)\ot \End_k(A)$ is a solution of the
Frobenius (resp. separability) equation.
\end{enumerate}
\end{theorem}

\begin{proof} $ 1)\Rightarrow 2)$
Both Frobenius and separable algebras are characterized by the existence
of an $A$-central element with some normalizing properties. Let
$e=\sum e^1\ot e^2\in A\ot A$ be such an $A$-central element. Then the map
$$R=R_e:\ A\ot A\to A\ot A, ~~~~R(a\ot b)=\sum e^1a\ot e^2 b$$
($a,b\in A$),
is a solution to the FS-equation. Here we view
$R_e\in \End_k(A\ot A)\cong \End_k(A)\ot \End_k(A)$ ($A$ is finitely
generated and projective over $k$).
Consequently, we can construct the algebra
${\cal A}(R)\subseteq \End_k(A)$.
We will prove that $A$ and ${\cal A}(R)$
are isomorphic when $A$ is a Frobenius algebra, or a separable algebra.\\
First we consider the injection $i:\ A\to \End_k(A)$, with $i(a)(b)=ab$, for
$a,b\in A$.
Then image of $i$ is included in ${\cal A}(R)$. Indeed,
$${\cal A}(R)=\{f\in \End_k(A)~|~(f\ot I_A)\circ R=R\circ (I_A\ot f)\}$$
Using the fact that $e$ is an $A$-central element, it follows easily that
$(i(a)\ot I_A)\circ R=R\circ (I_A\ot i(a))$, for all $a\in A$,
proving that ${\rm Im}(i)\subseteq {\cal A}(R)$.\\
If $f\in {\cal A}(R)$, then $(f\ot I_A)\circ R=R\circ (I_A\ot f)$, and,
evaluating this equality at $1_A\ot a$, we find
\begin{equation}\eqlabel{star}
\sum f(e^1)\ot e^2 a = \sum e^1 \ot e^2 f(a)
\end{equation}
Now assume that $A$ is a Frobenius algebra. Then there exists
$\eps:\ A\to k$ such that $\sum\eps(e^1)e^2=\sum e^1\eps(e^2)=1_A$.
Applying $\eps\ot I_A$ to \eqref{star} we obtain
that
$$f(a)=(\sum \eps(f(e^1))e^2)a$$
for all $a\in A$. Thus, $f=i(\sum \eps(f(e^1))e^2)$.
This proves that ${\rm Im}(i) = {\cal A}(R)$, and
the corestriction of $i$ to ${\cal A}(R)$ is an isomorphism of algebras.\\
If $A$ is separable, then $\sum e^1e^2=1_A$. Applying $m_A$
to \eqref{star} we find
$$f(a)=(\sum f(e^1)e^2)a$$
for all $a\in A$. Consequently $f=i(\sum f(e^1)e^2)$, proving again that
$A$ and ${\cal A}(R)$ are isomorphic.

$2) \Rightarrow 1)$ This is the last statement of \thref{4.1}.
\end{proof}

If $A$ is free of finite rank $n$ as a $k$-module, then we can describe
the algebra ${\cal A}(R)$ using generators and relations.
Let $M$ be free as a $k$-module, with basis
$\{m_1, m_2, \ldots, m_n\}$. In \seref{1.2}, we have seen that
an endomorphism
$R\in\End_k(M\ot M)$ is determined
by a matrix (see \eqref{matrix5} and \eqref{matrix6})
$R=(x^{ij}_{uv})_{i,j,u,v=1,n}\in {\cal M}_{n^2}(k)$.
Suppose that $R$ is a solution of
FS-equation. Identifying
$\End_k(M)$ and ${\cal M}_n(k)$, we will write ${\cal A}(n, R)$ for
the subalgebra of ${\cal M}_n(k)$ corresponding to ${\cal A}(R)$.
An easy computation shows that
\begin{equation}\eqlabel{3.2.1}
{\cal A}(n, R)=\{\left(a^i_{j}\right)\in {\cal M}_n(k)~|~
\sum_{\alpha=1}^n a^{\alpha}_ux^{ij}_{\alpha v}=
\sum_{\alpha=1}^n x^{i\alpha}_{u v}a^j_{\alpha}, \quad (\forall)\;
i, j, u, v=1,\cdots, n \}
\end{equation}
where $R$ is a matrix satisfying \eqref{3.3.1}.

\begin{corollary}\colabel{4.3}
Let $A$ be a $k$-algebra which is free of rank $n$ over $k$. Then
the following statements are equivalent:\\
1) $A$ is a Frobenius (resp. separable) algebra;\\
2) there exists an algebra isomorphism
$$A\cong {\cal A}(n, R),$$
where $R\in {\cal M}_{n^2}(k)\cong {\cal M}_{n}(k)\ot {\cal M}_{n}(k)$
is a solution of the Frobenius (resp. separability) equation.
\end{corollary}

\begin{remarks}\relabel{4.4}\rm
1) Over a field $k$ that is algebraically closed or of
characteristic zero, the structure of finite
dimensional separable $k$-algebras is given by the classical
Wedderburn-Artin Theorem: a finite dimensional algebra $A$
is separable if and only if is semisimple, if and only if
it is a direct product of matrix algebras.\\
2) As we have seen in \reref{3.13}, the dual $A^*$ of a separable finitely
generated projective algebra is a coseparable finitely generated projective
coalgebra. Thus we obtain a structure Theorem for coseparable
coalgebras, by using duality arguments.\\
More precisely,
let $C$ be a coseparable $k$-coalgebra which is free of rank $n$ over $k$.
Then there exists an
FS-map $\sigma:\ C\ot C\to k$ satisfying the normalizing condition
\eqref{3.10.2}. Let $A=C^*$ and $A\ot A\cong C^*\ot C^*$.
Then $\sigma$ is an $A$-central element of $A$, satisfying the
normalizing condition \eqref{SN}.
>From \coref{4.3}, we see that $A\cong {\cal A}(n,\sigma)$, and therefore
$C\cong {\cal A}(n,\sigma)^*$. Now ${\cal A}(n,\sigma)^*$
can be described as a quotient of the comatrix coalgebra:
 ${\cal A}(n,\sigma)^*={\cal M}^n(k)/I$, where $I$ is the coideal of
${\cal M}^n(k)$ that annihilates ${\cal A}(n,\sigma)$; $I$ is
generated by
$$\{o(i,j,u,v)=\sum_{\alpha=1}^n \left(c^{\alpha}_ux^{ij}_{\alpha v}-
x^{i\alpha}_{uv} c^j_{\alpha}\right)~~|~~i,j,k,l=1\ldots,m\}$$
Here we use the notation introduced in \seref{1.2}, and
$x^{ik}_{jl}=\sigma(c^i_j\ot c^k_l)$.
In a similar way, we can describe finite dimensional coalgebras $C$
for which the corestriction of scalars functor is Frobenius.
\end{remarks}

\begin{examples}\rm\exlabel{4.6}
1) Let $M$ be a free module of rank $n$, and $R=I_{M\ot M}$.
Then
$${\cal A}(R)=\{f\in \End_k(M)~|~f\ot I_M=I_M\ot f\}=k$$
2) Now let $R=\tau_M$ be the switch map. For all $f\in \End_k(M)$ and
$m,n\in M$, we have
$$((f\ot I_M)\circ\tau)(m\ot n)=f(n)\ot m=(\tau\circ (I_M\ot f))(m\ot n)$$
and, consequently, ${\cal A}(\tau)\cong {\cal M}_n(k)$ and
${\cal C}(\tau)\cong {\cal M}^n(k)$.\\
3)  Let $M$ be a finite dimensional vector space over a field $k$, and
$f\in \End_k(M)$ an idempotent. Then
$${\cal A}(f\ot f)=\{g\in \End_k(M)~|~f\circ g=g\circ f=\alpha f~
{\rm for~some~}\alpha\in k\}$$
Indeed,
$g\in {\cal A}(f\ot f)$ if and only  if $g\circ f\ot f=f\ot f\circ g$.
Multiplying the two factors, we find that $g\circ f=f\circ g$.
Now $g\circ f\ot f=f\ot g\circ f$ implies that $g\circ f=\alpha f$ for
some $\alpha\in k$. The converse is obvious.\\
In particular, assume that $M$ has dimension $2$ and let
$\{m_1,m_2\}$ be a basis of $M$. Let $f$ be the idempotent
endomorphism with matrix
$$\pmatrix{1-rq& q\cr r(1-rq)& rq\cr}$$
Assume first that $rq\not= 1$ and $r\not= 0$, and take $g\in \End_k(M)$
with matrix
$$\pmatrix{a& b\cr c& d}$$
$g\in {\cal A}(f\ot f)$ if and only if
$$
\alpha=a+br~~~;~~~
c+dr=r(a+br)~~~;~~~
br(1-rq)=qc$$
The two last equations can be easily solved for $b$ and $c$ in terms
of $a$ and $d$, and we see that ${\cal A}(f\ot f)$ has dimension two.
We know from the proof of \thref{4.1} that $f\in {\cal A}(f\ot f)$.
Another solution of the above system is
$I_M$, and we find that ${\cal A}(f\ot f)$ is the two-dimensional subalgebra
of $\End_k(M)$ with basis $\{f,I_M\}$. Put $f'=I_M-f$. Then $\{f,f'\}$
is also a basis for ${\cal A}(f\ot f)$, and ${\cal A}(f\ot f)=k\times k$.
${\cal C}(f\ot f)$ is the grouplike coalgebra of dimension two. We find the
same result if $rq=1$ or $q=0$.\\
4) Let $R\in {\cal M}_{n^2}(k)$ given by equation \eqref{3.14.4bis} as in
\prref{lazio}.
Then the algebra ${\cal A}(n, R)$ is given by
\begin{equation}
{\cal A}(n, R)=\{\left(a^i_{j}\right)\in {\cal M}_n(k)~|~
a^{\theta(j,v,i)}_u=
a^j_{\theta(u,i,v)}, \quad (\forall)\;
i, j, u, v\in X\}
\end{equation}
\prref{lazio} tells us when this algebra is separable
or Frobenius over $k$.\\
Assume now that $G$ is a finite group with $|G|=n$ invertible in $k$ and
$\theta$ is given as in \exref{3.15} 2).
Then the above algebra ${\cal A}(n, R)$ equals
$${\cal A}=\{\left(a^i_{j}\right)\in {\cal M}_n(k)~|~
a^i_{g(i)}=
a^j_{g(j)}, \quad (\forall)\;
i, j\in X~{\rm and}~g\in G \}$$
Indeed, if $(a^i_j)\in {\cal A}$, then $(a^i_j)\in {\cal A}(n,R)$, since
$g(\theta(j,v,i))=u$ implies that $g(j)=\theta(u,i,v)$. Conversely, for
$(a^i_j)\in {\cal A}(n,R)$, we choose $i$, $j\in X$ and $g\in G$.
Let $u=g(i)$ and $v=g(j)$. Then $\theta(j,v,u)=i$, hence
$$
a^i_{g(i)}=a^{\theta(j,v,u)}_u=a^j_{\theta(u,u,v)}=a^j_v=a^j_{g(j)},
$$
showing that ${\cal A}={\cal A}(n,R)$.
>From the fact that $G$ acts transitively,
it follows that a matrix in ${\cal A}(n,R)$ is completely determined
by its top row. For every $g\in G$, we define $A_g\in {\cal A}(n,R)$
by $(A_g)^1_i=\delta_{g(1),i}$. Then ${\cal A}(n,R)=\{A_g~|~g\in G\}$,
and we have an algebra isomorphism
$$f:\ {\cal A}(n, R) \to kG, \quad f(A_g)= g$$
For example, take the cyclic group of order $n$, $G=C_n$.
$${\cal A}(n, R)=\left\{\pmatrix{
 x_1& x_2& \cdots &x_n \cr
 x_n& x_1& \cdots &x_{n-1} \cr
 x_{n-1}& x_n& \cdots & x_{n-2}\cr
 \cdot & \cdot&\cdots &\cdot\cr
x_2& x_3&\cdots & x_1}~|~
x_1, \cdots, x_n\in k\right\}\cong kC_n$$
5) Let $G$ be a finite group of order $n$ and $\sigma:\ G\times G\to k^*$
a $2$-cocycle. Let $R_{\sigma}$ be the solution of the
$FS$-equation given by \eqref{cocsol}. We then obtain directly from
\eqref{3.2.1} that ${\cal A}(n, R_{\sigma})$ consists
of all $G\times G$-matrices $(a_j^i)_{i,j\in G}$ satisfying the relations
$$
a_u^{jv^{-1}i}\; \sigma ^{-1} (ui^{-1}vj^{-1}, jv^{-1})
\sigma (vj^{-1}, jv^{-1}i) \sigma (jv^{-1}, v)=
$$
$$
a_{ui^{-1}v}^j\; \sigma ^{-1} (iu^{-1}, ui^{-1})
\sigma (iu^{-1}, u) \sigma (ui^{-1}, v)
$$
for all $i$, $j$, $u$, $v\in G$. This algebra is separable
if $n$ is invertible in $k$.
\end{examples}

We will now present some new classes of examples, starting from the
solution $R^{\phi}$ of the $FS$-equation discussed in \prref{3.16}.
In this case, we find easily that
$$x^{ij}_{uv}=\delta_{uv}\delta_{\phi(i)u}\delta_{\phi(j)v}$$
and, according to \eqref{3.2.1}, ${\cal A}(R^{\phi})$ consists of matrices
$\left(a^i_j\right)$ satisfying
$$\sum_{\alpha=1}^n a^{\alpha}_u x^{ij}_{\alpha v}=
\sum_{\alpha=1}^n  x^{i\alpha}_{u v}a^j_{\alpha}$$
or
\begin{equation}\eqlabel{3.4.1}
a^v_u\delta_{\phi(i)v}\delta_{\phi(j)v}= \sum_{\alpha\in\phi^{-1}(v)}
\delta_{uv} \delta_{\phi(i)u} a_{\alpha}^j
\end{equation}
for all $i,j,v,u=1,\ldots,n$. The left hand side of \eqref{3.4.1} is nonzero
if and only if $\phi(i)=\phi(j)=v$. If $\phi(i)=\phi(j)=v=u$, then
\eqref{3.4.1}
 amounts to
$$a_{\phi(i)}^{\phi(i)}=\sum_{\{\alpha|\phi(\alpha)=\phi(i)\}}
a^i_{\alpha}.$$
If $\phi(i)=\phi(j)=v\neq u$, then \eqref{3.4.1} takes the form
$$a^{\phi(i)}_u=0.$$
Now assume that the left hand side of \eqref{3.4.1} is zero. If
$\phi(i)\neq \phi(j)$,
then the right hand side of \eqref{3.4.1} is also zero, except when
$u=v=\phi(i)$. Then \eqref{3.4.1} yields
$$\sum_{\{\alpha|\phi(\alpha)=\phi(i)\}} a^j_{\alpha}=0.$$
If $\phi(i)=\phi(j)\neq u$, then \eqref{3.4.1} reduces to $0=0$.
We summarize our results as follows.

\begin{proposition}\prlabel{4.7}
Consider an idempotent map $\phi:\ \{1,\ldots,n\}\to  \{1,\ldots,n\}$,
and the corresponding solution $R^{\phi}$ of the $FS$-equation. Then
${\cal A}(R^{\phi})$ is the subalgebra of ${\cal M}_n(k)$ consisting of
matrices satisfying
\begin{eqnarray}
a_{\phi(i)}^{\phi(i)}&=&\sum_{\{\alpha|\phi(\alpha)=\phi(i)\}} a^i_{\alpha}
~~~~~~~~~~~~(i=1,\ldots,n)\eqlabel{3.4.2}\\
a^{\phi(i)}_j&=&\hspace*{1cm}0\hspace*{28mm}(\phi(i)\neq j)
\eqlabel{3.4.3}\\
0&=&\sum_{\{\alpha|\phi(\alpha)=\phi(i)\}} a^j_{\alpha}
~~~~~~~~~~~~(\phi(i)\neq \phi(j))\eqlabel{3.4.4}
\end{eqnarray}
${\cal A}(R^{\phi})$ is a separable
$k$-algebra if and only if $(R^{\phi},\varepsilon={\rm trace})$
is a solution of the F-equation
if and only if $\phi$ is the identity map. In this case,
${\cal A}(R^{\phi})$ is the direct sum of $n$ copies of $k$.
\end{proposition}

\begin{proof}
The first part was done above. ${\cal A}(R)$ is separable if and only if
\eqref{3.3.2} holds. This comes down to
$$\delta_{ju}=\sum_{v=1}^n \delta_{uv}\delta_{\phi(u)u}\delta_{\phi(j)v}
=\delta_{\phi(u)u}\delta_{\phi(j)u}$$
and this implies that $\phi(u)=u$ for all $u$.
(\ref{eq:3.4.2},\ref{eq:3.4.3},\ref{eq:3.4.4}) reduce to $a^i_j=0$ for
$i\not= j$, and ${\cal A}(R^I)$ consists of all diagonal matrices.\\
$(R^{\phi},\varepsilon={\rm trace})$
is a solution of the F-equation if and only if \eqref{3.3.3} holds,
and a similar computation shows that this also implies that $\phi$ is
the identity.
\end{proof}

\begin{examples}\rm\exlabel{4.8}
1) Take $n=4$, and $\phi$ given by
$$\phi(1)=\phi(2)=2,~\phi(3)=\phi(4)=4$$
(\ref{eq:3.4.2},\ref{eq:3.4.3},\ref{eq:3.4.4}) take the following form
$$\begin{array}{ccc}
a_2^2=a_1^1+a_2^1&&a_4^4=a_3^3+a_4^3\\
a_1^2=a_2^2=a_4^2=0&&a_1^4=a_2^4=a_3^4=0\\
a_1^3=-a_2^3&&a_3^1=-a_4^1
\end{array}$$
and
$${\cal A}(4, R^{\phi})=\left\{
\pmatrix{x& y-x& u&-u\cr
0&y&0&0\cr
v&-v&z&t-z\cr
0&0&0&t\cr}~|~
x, y, z, t, u, v\in k\right\}$$
The dual coalgebra can also be described easily. Write
$x_i=c_i^i$ ($i=1,\ldots,4$), $x_5=c_3^1$ and $x_6=c_1^3$. Then
${\cal C}(R^{\phi})$ is the six dimensional coalgebra with basis
$\{x_1,\ldots,x_6\}$ and
\begin{eqnarray*}
\Delta(x_1)&=& x_1\ot x_1+ x_5\ot x_6, \quad
\Delta(x_2)= x_2\ot x_2,\\
\Delta(x_3)&=& x_3\ot x_3+ x_6\ot x_5, \quad
\Delta(x_4)= x_4\ot x_4,\\
\Delta(x_5)&=& x_1\ot x_5+ x_5\ot x_3, \quad
\Delta(x_6)= x_6\ot x_1+ x_3\ot x_6,\\
\varepsilon(x_i)&=& 1~~~~~(i=1,2,3,4), \quad
\varepsilon(x_i)=0~~~~~(i=5,6).
\end{eqnarray*}
2) Again, take $n=4$, but let $\phi$ be given by the formula
$$\phi(1)=1,~\phi(2)=\phi(3)=\phi(4)=2$$
Then (\ref{eq:3.4.2},\ref{eq:3.4.3},\ref{eq:3.4.4}) reduce to
$$a_2^2=a_2^3+a_3^3+a_4^3=a_2^4+a_3^4+a_4^4$$
$$a_2^1=a_3^1=a_4^1=a_1^2=a_3^2=a_4^2=a_1^3=a_1^4=0,$$
hence
$${\cal A}(4, R^{\phi})=\left\{
\pmatrix{x&0&0&0\cr
0&y&0&0\cr
0& u&z& y-z-u\cr
0&v&y-t-v& t\cr}
~|~
x, y, z, t, u, v\in k\right\}$$
Putting $c_i^i=x_i$ ($i=1,\ldots,4$), $x_5=c_2^3$ and $x_6=c_2^4$, we find
that ${\cal C}(R^{\phi})$ is the six dimensional coalgebra with structure
maps
\begin{eqnarray*}
\Delta(x_1)&=& x_1\ot x_1, \quad
\Delta(x_2)= x_2\ot x_2,\\
\Delta(x_3)&=& x_3\ot x_3+ (x_2-x_3-x_5)\ot (x_2-x_4-x_6),\\
\Delta(x_4)&=& x_4\ot x_4+ (x_2-x_4-x_6)\ot (x_2-x_3-x_5),\\
\Delta(x_5)&=& x_5\ot x_2+ x_3\ot x_5+ (x_2-x_3-x_5) \ot x_6,\\
\Delta(x_6)&=& x_6\ot x_2+ x_4\ot x_6+(x_2-x_4-x_6)\ot x_5,\\
\varepsilon(x_i)&=& 1~~~~~(i=1,2,3,4), \quad
\varepsilon(x_i)=0~~~~~(i=5,6).
\end{eqnarray*}
3) Put $\phi(i)=1$, for all $i=1,\ldots,n$.
(\ref{eq:3.4.2},\ref{eq:3.4.3},\ref{eq:3.4.4}) reduce to
$$a^1_i=0~~~{\rm and}~~ \sum_{\alpha=1}^n a^i_{\alpha}=a^1_1,$$
for all $i\neq 1$.
\end{examples}

\section{The category of FS-objects}\selabel{5}
Our starting point is the following result, due to
Abrams (cf.  \cite[Theorem 3.3]{A3}):
over a finite dimensional Frobenius algebra $A$ the categories
${\cal M}_A$ and ${\cal M}^A$ are isomorphic. In this Section,
we will extend this result to weak Frobenius algebras.\\
We have seen in \coref{3.4} that the equation
$R^{12}R^{23}=R^{13}R^{12}$ is equivalent to the fact that a certain
multiplication on $M\ot M$ is associative. We shall now prove
that the other equation, namely $R^{12}R^{23}=R^{23}R^{13}$, is equivalent
to the fact that a certain comultiplication is coassociative.

\begin{proposition}\thlabel{5.1}
Let $(A, m_A, 1_A)$ be an algebra, $R=\sum R^1\ot R^2\in A\ot A$ and
$$
\delta :A\to A\ot A,\quad \delta(a)=\sum R^1\ot R^2a
$$
for all $a\in A$. The following statements are equivalent
\begin{enumerate}
\item[1)] $(A,\delta)$ is a coassociative coalgebra (not necessarily
with a counit).
\item[2)] $R^{12}R^{23}=R^{23}R^{13}$ in $A\ot A\ot A$.
\end{enumerate}
In this case any left $A$-module $(M,\cdot)$ has a structure of left
comodule over the coalgebra $(A,\delta)$ via
$$\rho:\ M\to A\ot M, \quad \rho(m):=\sum R^1\ot R^2\cdot m$$
for all $m\in M$.
\end{proposition}

\begin{proof}
The equivalence of 1) and 2) follows from the formulas
\begin{eqnarray*}
(\delta\ot I)\delta(a)&=&R^{12}R^{23}\cdot (1_A\ot 1_A\ot a)\\
(I\ot \delta)\delta(a)&=&R^{23}R^{13}\cdot (1_A\ot 1_A\ot a)
\end{eqnarray*}
for all $a\in A$. The final statement follows from
\begin{eqnarray*}
(\delta\ot I)\rho(m)&=&R^{12}R^{23}\cdot (1_A\ot 1_A\ot m)\\
(I\ot \rho)\rho(m)&=&R^{23}R^{13}\cdot (1_A\ot 1_A\ot m)
\end{eqnarray*}
for all $m\in M$.
\end{proof}

Suppose now that $(A, m_A, 1_A)$ is an algebra over $k$ and let $R\in A\ot A$
be a $A$-central element. Then $R^{12}R^{23}=R^{23}R^{13}$ in
$A\ot A\ot A$. It follows that
$(A,\Delta_R:=\delta^{{\rm cop}})$ is also a coassociative coalgebra,
where $\Delta_R(a)=\delta^{{\rm cop}}(a)=\sum R^2a\ot R^1$, for all
$a\in A$. We remark that $\Delta_R$ is not an algebra map, i.e.
$(A,m_A,\Delta_R)$ is not a bialgebra. Any left $A$-module $(M,\cdot)$ has a
structure of right comodule over the coalgebra $(A,\Delta_R)$ via
$$
\rho_R:\ M\to M\ot A, \quad \rho_R(m):=\sum R^2\cdot m\ot R^1
$$
for all $m\in M$.
Moreover, for any $a\in A$ and $m\in M$ we have that
$$
\rho_R(a\cdot m)=\sum a_{(1)}\cdot m\ot a_{(2)}=
\sum m_{<0>}\ot am_{<1>}.
$$
Indeed, from the definition of the coaction on $M$ and the comultiplication
on $A$, we have immediately
$$
\rho_R(a\cdot m)=\sum R^2a\cdot m\ot R^1=\sum a_{(1)}\cdot m\ot a_{(2)}$$
On the other hand
$$
\sum m_{<0>}\ot am_{<1>}=\sum R^2\cdot m\ot aR^1=\sum R^2a\cdot m\ot R^1
$$
where in the last equality we used that $R$ is a $A$-central element. These
considerations lead us to the following

\begin{definition}\delabel{5.2}
Let $(A, m_A, \Delta_A)$ be at once an algebra and a coalgebra
(but not necessarily a bialgebra).
An $FS$-object over $A$ is a $k$-module $M$ that is at once a left
$A$-module and a right $A$-comodule such that
\begin{equation}\eqlabel{5.2.1}
\rho(a\cdot m)=\sum a_{(1)}\cdot m\ot a_{(2)}=\sum m_{<0>}\ot am_{<1>}
\end{equation}
for all $a\in A$ and $m\in M$.
\end{definition}

The category of $FS$-objects and $A$-linear $A$-colinear maps will be
denoted by ${}_A{\cal FS}^A$. This category measures how far
$A$ is  from a bialgebra. If $A$   has not
a unit (resp. a counit), then the objects in ${}_A{\cal FS}^A$ will be
assumed to be unital (resp. counital).

\begin{proposition}\prlabel{5.3}
If $(A, m_A, 1_A, \Delta_A, \varepsilon_A)$ is a bialgebra with
unit and counit, then the forgetful functor
$$F:{}_A{\cal FS}^A \to {}_k{\cal M}$$
is an isomorphism of categories.
\end{proposition}

\begin{proof}
Define $G:\ {}_k{\cal M}\to {}_A{\cal FS}^A$ as follows: $G(M)=M$ as a
$k$-module, with trivial $A$-action and $A$-coaction:
$$\rho(m)=m\ot 1_A~~{\rm and}~~a\cdot m=\varepsilon_A(a)m$$
for all $a\in A$ and $m\in M$. It is clear that
$G(M)\in {}_A{\cal FS}^A$.\\
Now, assume that $M$ is an $FS$-object over $A$. Applying
$I\ot\varepsilon_A$
to \eqref{5.2.1}, we find that
$$a\cdot m=\sum m_{<0>}\varepsilon_A(a)\varepsilon_A(m_{<1>})=
\varepsilon_A(a)m$$
Taking $a=1_A$ in \eqref{5.2.1}, we find
$$\rho(m)=\sum 1_A\cdot m\ot 1_A=m\ot 1_A.$$
Hence, $G$ is an inverse for the forgetful functor
$F:{}_A{\cal FS}^A \to {}_k{\cal M}$.
\end{proof}

\begin{definition}\delabel{5.4}
A triple $(A, m_A, \Delta_A)$ is called a weak Frobenius algebra
($WF$-algebra, for short) if \\ 
$(A, m_A)$ is an algebra
(not necessarily with unit), $(A, \Delta_A)$ is a coalgebra
(not necessarily with counit) and\\
$(A, m_A, \Delta_A)\in {}_A{\cal FS}^A$, that is
\begin{equation}\eqlabel{5.3.1}
\Delta(ab)=\sum a_{(1)}b\ot a_{(2)}=\sum b_{(1)}\ot ab_{(2)}.
\end{equation}
for all $a$, $b\in A$.
\end{definition}

\begin{remarks}\relabel{dumi}\relabel{5.5}
\rm
1) Assume that $A$ is an $WF$-algebra with unit, and write
$\Delta(1_A)=\sum e^2\ot e^1$. From \eqref{5.3.1}, it follows that
\begin{equation}\eqlabel{4.1.3}
\Delta(a)=\sum e^2a\ot e^1=\sum e^2\ot ae^1
\end{equation}
and this implies that $\Delta^{\rm cop}(1_A)=\sum e^1\ot e^2$ is an
$A$-central element.
Conversely, if $A$ is an algebra with unit, and $e=\sum e^1\ot e^2$
is an $A$-central element, then $e^{12}e^{23}=e^{23}e^{13}$ (see
\eqref{2.3.2a}),
and it is easy to prove that this last statement is equivalent to the
fact that $\Delta:\ A\to A\ot A$ given by $\Delta(a)=\sum e^2a\ot e^1$
is coassociative. Thus $A$ is a $WF$-algebra.
We have proved that $WF$-algebras with unit correspond to
algebras with unit together with an $A$-central element.\\
2) From \eqref{5.3.1}, it follows immediately that
$f:=\Delta^{\rm cop}:\ A\to A\ot A$ is an $A$-bimodule map.
Conversely, if $f$ is an $A$-bimodule map, then it is easy to prove
that $\Delta=\tau\circ f$ defines a coassociative comultiplication on
$A$, making $A$ into a $WF$-algebra. Now, using \coref{adi}, we obtain
that a finitely generated projective and unitary $k$-algebra
$(A, m_A, 1_A)$ is Frobenius if and only if $A$
is an unitary and counitary $WF$-algebra. Thus, we can view $WF$-algebras
as a generalization of Frobenius algebras.
\end{remarks}

\begin{proposition}\prlabel{5.6}
Let $(A, m_A, 1_A, \Delta_A)$ be a $WF$-algebra with unit. Then
the forgetful functor
$$F:\ {}_A{\cal FS}^A \to {}_A{\cal M}$$
is an
isomorphism of categories.
\end{proposition}

\begin{proof}
We define a functor $G:\ {}_A{\cal M}\to{}_A{\cal FS}^A$ as follows:
$G(M)=M$ as an $A$-module, with right $A$-coaction given by the formula
$$\rho(m)=\sum e^2\cdot m\ot e^1$$
where $\Delta(1_A)=\sum e^2\ot e^1$. $\rho$ is a coaction because
$e^{12}e^{23}=e^{23}e^{13}$, and, using \eqref{4.1.3}, we see that
\begin{eqnarray*}
\rho(a\cdot m)&=& \sum e^2a\cdot m\ot e^1=\sum a_{(1)}\cdot m\ot a_{(2)}\\
&=& \sum e^2\cdot m\ot ae^1=\sum m_{<0>}\ot am_{<1>}
\end{eqnarray*}
as needed. Now, $G$ and $F$ are each others inverses.
\end{proof}

Now, we will give the coalgebra version of \prref{5.6}.
Consider a $WF$-algebra $(A, m_A, \Delta_A, \varepsilon_A)$, with a counit
$\varepsilon_A$, and consider
$\sigma=\varepsilon\circ m_A\circ \tau:\ A\ot A\to k$, that is,
$$\sigma(c\ot d)= \varepsilon(dc)$$
for all $c,d\in A$. Now
$$\Delta(cd)=\sum c_{(1)}d\ot c_{(2)}=\sum d_{(1)}\ot cd_{(2)}$$
so
$$\sum \sigma(d\ot c_{(1)})c_{(2)}=
(\varepsilon\ot I_C)(\Delta(cd))=
 (I_C\ot\varepsilon)(\Delta(cd))=\sum \sigma(d_{(2)}\ot c)d_{(1)}$$
and $\sigma$ is an $FS$-map. Conversely, let $(A, \Delta_A)$ be a
coalgebra with counit, and assume that $\sigma:\ A\ot A\to k$
is an $FS$-map. A straightforward computation shows that the formula
$$c\cdot d=\sum \sigma(d_{(2)}\ot c)d_{(1)}$$
defines an associative multiplication on $A$ and that
$(A, \cdot, \Delta_A)$ is a $WF$-algebra. Thus,
$WF$-algebras with counit correspond to coalgebras with counit
together with an $FS$-map.

\begin{proposition}\prlabel{5.7}
Let $(A, m_A, \Delta_A, \varepsilon_A)$ be a $WF$-algebra with counit.
Then  the forgetful functor
$$F:{}_A{\cal FS}^A \to {\cal M}^A$$
is an isomorphism of categories.
\end{proposition}

\begin{proof} The inverse of $F$ is the
functor $G:\ {\cal M}^A\to {}_A{\cal FS}^A$ defined as follows:
$G(M)=M$ as a $A$-comodule, with $A$-action given by
$$a\cdot m=\sum \sigma(m_{<1>}\ot a)m_{<0>}$$
for all $a\in A$, $m\in M$. Further details are left to the reader.
\end{proof}

As an immediate consequence of \prref{5.6} and \prref{5.7}
we obtain the following generalization of Abrams' result
\cite[Theorem 3.3]{A3}

\begin{corollary}\colabel{5.8}
Let $(A, m_A, 1_A, \Delta_A, \varepsilon_A)$ be a
$WF$-algebra with unit and counit. Then we have an equivalence of
categories
$${}_A{\cal FS}^A\cong {}_A{\cal M} \cong {\cal M}^A$$
\end{corollary}

Let us finally show that \prref{5.6} also holds over
$WF$-algebras that are unital as modules over themselves.
Examples of such algebras are the central separable algebras
discussed in \cite{T}.

\begin{proposition}\prlabel{5.9}
Let $A$ be a $WF$-alegebra that is unital as a module over itself. We have
an equivalence between the categories ${}_A{\cal M}$ and ${}_A{\cal FS}^A$.
\end{proposition}

\begin{proof}
For a unital $A$-module $M$, we define $F(M)$ as the $A$-module $M$ with
$A$-coaction given by
$$\rho(a\cdot m)= \sum a_{(1)}\cdot m\ot a_{(2)}$$
It is clear that $\rho$ defines an $A$-coaction. One equality in
\eqref{5.2.1} is obvious, and the other one follows from \eqref{5.3.1}:
for all $a,b\in A$ and $m\in M$, we have that
\begin{eqnarray*}
\sum (b\cdot m)_{(0)}\ot a(b\cdot m)_{(1)}&=&
\sum b_{(1)}\cdot m\ot ab_{(2)}\\
&=& \sum (ab)_{(1)}\cdot m\ot (ab)_{(2)}\\
&=& \rho(a\cdot(b\cdot m))
\end{eqnarray*}
It follows that $F(M)$ is an $FS$-object, and $F$ defines the desired
category equivalence.
\end{proof}

{\bf Acknowledgement}
The second author thanks Lowell Abrams for fruitful discussions.
We also are indebted to Bogdan Ichim for writing the computer program
for the numerical computations in \exref{3.5}, and to Mona St\u{a}nciulescu
for her stimulating comments to the early versions of this paper.


\begin{thebibliography}{10}

\bibitem{Ab}
E. Abe, Hopf algebras, Cambridge Univ. Press, Cambridge, 1977.

\bibitem{A1}
L. Abrams, Two-dimensional topological quantum field theories
and Frobenius algebras, {\sl J. Knot Theory Ramifications} {\bf 5}
(1996), 569-587.

\bibitem{A2}
L. Abrams, The quantum Euler class and the quantum cohomology of the
Grassmanians, preprint q-alg/9712025.

\bibitem{A3}
L. Abrams, Modules, comodules and cotensor products over Frobenius
algebras, {\sl J. Algebra}, to appear.

\bibitem{A}
M. Atiyah, An introduction to topological quantum field theories,
{\sl Turkish J. Math.} {\bf 21} (1997), 1-7.

\bibitem{Ba1}
R. J. Baxter, Exactly solved models in statistical mechanics,
Academic Press, London, 1989.

\bibitem{BDGN}
M. Beattie, S. D\u asc\u alescu, L. Gr\"unenfelder and
 C. N\u ast\u asescu,
Finiteness conditions, co-Frobenius Hopf algebras and quantum
groups, {\sl J. Algebra} {\bf 200} (1998), 312-333.

\bibitem{BFS0}
K.I. Beidar, Y. Fong and A. Stolin, On antipodes and integrals in Hopf
algebras over rings and the quantum Yang-Baxter equation,
{\sl J. Algebra} {\bf 194} (1997), 36-52.

\bibitem{BFS}
K.I. Beidar, Y. Fong and A. Stolin, On Frobenius algebras and the
Yang-Baxter equation, {\sl Trans. Amer. Math. Soc.} {\bf 349} (1997),
3823-3836.

\bibitem{BF}
A.D. Bell and R. Farnsteiner, On the theory of Frobenius extensions and its
application to Lie superalgebras, {\sl Trans. Amer. Math. Soc.} {\bf 335}
(1993),
407-424.

\bibitem{B1}
T. Brzezi\'nski, Frobenius properties and Maschke-type theorems
for entwined modules, {\sl Proc. Amer. Math. Soc.}, to appear.

\bibitem{B2}
T. Brzezi\'nski, Coalgebra-Galois extensions from the extension theory
point of view, to appear in ``Hopf algebras and quantum groups",
S. Caenepeel and F. Van Oystaeyen (Eds.), {\sl Lecture Notes Pure
Appl. Math.}, Dekker, New York, 2000.

\bibitem{BM}
T. Brzezi\'nski, S. Majid, Coalgebra bundles, {\sl Comm. Math. Phys.},
{\bf 191} (1998), 467-492.

\bibitem{CMZ}
S. Caenepeel, G. Militaru and S. Zhu, Doi-Hopf modules, Yetter-Drinfel'd
modules and Frobenius type properties, {\sl Trans. Amer. Math. Soc.}
{\bf 349} (1997), 4311-4342.

\bibitem{CIMZ}
S. Caenepeel, Bogdan Ion, G. Militaru and S. Zhu,
Separable functors for Doi-Hopf modules. Applications, {\sl Advances Math.}
{\bf 145} (1999), 239-290.

\bibitem{CGN}
F. Casta\~{n}o Iglesias, J. G\'omez Torrecillas and
C. N\u ast\u asescu, Frobenius Functors. Applications, {\sl Comm. Algebra},
to appear.

\bibitem{CR}
C. Curtis and I. Reiner, Representation theory of finite groups and
associative algebras, {\sl Interscience Publishers}, New York, 1962.

\bibitem{D}
Y. Doi, Homological coalgebra, {\sl J. Math. Soc. Japan}
{\bf 33} (1981), 31-50.

\bibitem{DI}
F. DeMeyer and E. Ingraham, Separable algebras over commutative rings,
{\sl Lecture Notes in Math.} {\bf 181}, Springer Verlag, Berlin, 1971.

\bibitem{Dr}
V. G. Drinfel'd, Hamiltonian structures on Lie groups, Lie bialgebras,
and the geometrical meaning of classical Yang-Baxter equation,
{\sl Soviet Math. Dokl.} {\bf 27} (1983), 68-71.

\bibitem{Dub}
B. Dubrovin, Geometry of 2d topological field theory,
{\sl Lecture Notes in Math.} {\bf 1620}, Springer Verlag, Berlin, 1996,
120-348.

\bibitem{ESS}
P. Etingof, T. Schedler, A. Soloviev, Set-theoretical solutions to the
quantum Yang-Baxter equation, {\sl Duke Math. J.}, to appear.

\bibitem{Fa}
L.D. Faddeev, Quantum completely integrable models in quantum field
theory, {\sl Sov. Sci. Rev}, {\bf C1} (1980), 107-155.

\bibitem{F}
R. Farnsteiner, On Frobenius extensions defined by Hopf algebras,
{\sl J. Algebra} {\bf 166} (1994), 130-141.

\bibitem{FMS}
D. Fischman, S. Montgomery and  H.J. Schneider, Frobenius extensions of
subalgebras of Hopf algebras, {\sl Trans. Amer. Math. Soc.} {\bf 349} (1997),
4857-4895.

\bibitem{Hi}
J. Hietarinta, Permutation-type solutions to the quantum Yang-Baxter
and other n-simplex equations, q-alg/9702006, to appear in
{\sl J. Phys. A: Math. Gen.}.

\bibitem{HS}
K. Hirata, K. Sugano, On semisimple and separable extensions over
noncommutative rings, {\sl J. Math. Soc. Japan} {\bf 18} (1966),
360-373.

\bibitem{JS}
A. Joyal, R. Street, An introduction to Tannaka duality and quantum
groups, {\sl Lecture Notes Math.} {\bf 1488} (1991), 411-492.

\bibitem{K1}
L. Kadison, The Jones polynomial and certain separable Frobenius
extensions, {\sl J. Algebra} {\bf 186} (1996), 461-475.

\bibitem{Kan}
D. M. Kan, Adjoint functors, {\sl Trans. Amer. Math. Soc.}
{\bf 87} (1958), 294-323.

\bibitem{FK}
F. Kasch, Grundlagen einer theory der Frobeniuserweiterungen,
{\sl Math. Ann.} {\bf 127} (1954), 453-474.

\bibitem{K} C. Kassel, Quantum Groups, {\sl Graduate Texts in Math.}
{\bf 155}, Springer Verlag, Berlin, 1995.

\bibitem{LR}
L.A. Lambe, D. Radford, Algebraic aspects of the
quantum Yang-Baxter equation, {\sl J. Algebra} {\bf 54} (1992), 228-288.

\bibitem{L}
R.G. Larson, Coseparable coalgebras, {\sl J. Pure and Appl. Algebra}
{\bf 3} (1973), 261-267.

\bibitem{Lin}
B.J. Lin, Semiperfect coalgebras, {\sl J. Algebra} {\bf 49} (1977),
357-373.

\bibitem{Maj 1}
S. Majid, Foundation of quantum group theory, Cambridge Univ. Press, 1995.

\bibitem{MN}
C. Menini and C. N\u ast\u asescu, When are induction and coinduction
functors isomophic?, {\sl Bull. Belg. Math. Soc.- Simon Stevin} {\bf 1} (1994),
521-558.

\bibitem{M3}
G. Militaru, The Hopf modules category and the pentagonal equation,
{\sl Comm. Algebra} {\bf 10} (1998), 3071-3097.

\bibitem{M4}
G. Militaru, A class of solutions for the integrability condition of the
Knizhnik-Zamolodchikov equation: a Hopf algebraic approach,
{\sl Comm. Algebra} {\bf 27} (1999), 2393-2407.

\bibitem{Montgomery}
S. Montgomery, Hopf algebras and their actions on rings, American Mathematical
Society, Providence, 1993.

\bibitem{Mo}
K. Morita, Adjoint pairs of functors and Frobenius extensions,
{\sl Sci. Rep. Tokyo Kyoiku Daigaku, Sect. A9} (1965), 40-71.

\bibitem{N1}
T. Nakayama and T. Tsuzuku, On Frobenius extensions I, {\sl Nagoya Math. J.}
{\bf 17} (1960), 89-110.

\bibitem{N2}
T. Nakayama and T. Tsuzuku, On Frobenius extensions II, {\sl Nagoya Math. J.}
{\bf 19} (1961), 127-148.

\bibitem{NOB}
C. N\u ast\u asescu, M. van den Bergh and F. van Oystaeyen, Separable functors
applied to graded rings, {\sl J. Algebra} {\bf 123} (1989), 397-413.

\bibitem{P}
B. Pareigis, When Hopf algebras are Frobenius algebras, {\sl J. Algebra}
{\bf 18} (1971), 588-596.

\bibitem{PAR}
B. Pareigis, Einige Bemerkungen uber Frobenius-Erweiterungen,
{\sl Math. Ann.} {\bf 153} (1964), 1-13.

\bibitem{FRT}
N. Yu Reshetikhin, L.A. Takhtadjian and L.D. Faddeev,
Quantization of Lie groups
and Lie algebras, {\sl Algebra i Analiz} {\bf 1} (1989), 178-206.

\bibitem{T}
J.L. Taylor, A bigger Brauer group, {\sl Pacific J. Math.} {\bf 103}
(1982), 163-203.

\bibitem{W}
Y. Watatani, Index of $C^*$-subalgebras, {\sl Memoirs Amer. Math. Soc.}
{\bf 83} (1990).

\end{thebibliography}
\end{document}